\documentclass[preprint,12pt]{elsarticle}
\biboptions{sort&compress}
\usepackage{amssymb}
\usepackage{amsmath}
\usepackage{algorithmic}
\usepackage{algorithm}
\usepackage{array}
\usepackage[caption=false,font=normalsize,labelfont=sf,textfont=sf]{subfig}
\usepackage{textcomp}
\usepackage{stfloats}
\usepackage{url}
\usepackage{verbatim}
\usepackage{graphicx}
\usepackage{xcolor}
\usepackage{rotating}
\usepackage{pgfplots}
\usepackage{tikz}
\pgfplotsset{compat=1.18}
\usepgfplotslibrary{groupplots}
\usepackage{booktabs}

\newcommand{\rot}[1]{\begin{turn}{60}#1\enspace\end{turn}}

\journal{Renewable Energy Focus}

\begin{document}

\begin{frontmatter}

\title{Transmission Expansion and Energy Storage Planning under Great Britain's Uncertain Energy Transition Pathways}

\author[strath]{Shengming Cao}
\author[strath]{Runsheng He}
\author[strath]{Waqquas Bukhsh\corref{cor1}}
\ead{waqquas.bukhsh@strath.ac.uk}
\cortext[cor1]{Corresponding author.}
\address[strath]{Department of Electronic and Electrical Engineering,
University of Strathclyde, Glasgow G1 1XQ, United Kingdom}

\begin{abstract}
Long term power system planning must account for uncertainty in electricity demand and generation, while recognising the different investment and commissioning characteristics of transmission and energy storage. This paper develops a multi year stochastic framework for coordinated transmission expansion and energy storage planning and applies it to the Great Britain electricity system to assess investment needs through to 2050. The framework integrates multiple energy transition pathways with an enhanced representation of the GB transmission network and captures the different characteristics of transmission and storage investment over time. The coordinated treatment of transmission and storage is considered as a system level benchmark rather than a centralised investment mechanism. In GB, transmission investment is subject to regulated forward looking planning, while storage is largely developed by market participants in response to market signals. The framework therefore provides a means of quantifying the potential value of improved coordination and identifying where existing investment signals may not encourage storage deployment that complements network development. Results across alternative transmission only, storage only and coordinated investment strategies show substantial benefits from considering the two asset types jointly. Compared with transmission only planning, coordinated investment reduces additional transmission capacity by 2\%, cumulative expected load shedding by 97\%, and total expected system cost by 36\%, primarily through avoided unserved energy. The results also show that transmission preparation periods affect the timing, but not the final portfolio, of investments. Finally, the energy transition pathway materially influences investment requirements, with Electrification and Efficiency requiring 7\% more transmission capacity and 45\% more storage energy capacity than Falling Behind. The framework provides a system level basis for understanding how network and storage investment can interact under uncertain energy transitions.
\end{abstract}

%

\end{frontmatter}

\section{Introduction}\label{sec:intro}

The transition to a low carbon energy system is central to meeting climate targets, improving energy security and supporting sustainable economic development~\cite{kabeyi2022sustainable}. Electricity networks are particularly important because the decarbonisation of transport, heating and industry will depend increasingly on electrification~\cite{bogdanov2021full}. This transition will substantially change where electricity is generated and consumed, creating new requirements for network capacity, flexibility and investment. Long term power system planning must therefore account for both the scale and uncertainty of future energy system development~\cite{shahzad2024renewable,cao2025stochastic,NESO_FES2025}.

The Great Britain (GB) electricity system is undergoing a major transformation towards the UK’s 2050 net zero target. This includes rapid growth in renewable generation and increasing electrification of transport, heating and industrial demand~\cite{dixon2022way,shahzad2024charting,royapoor2023towards}. These changes are expected to have a major impact on the transmission network over the coming decades. In particular, the period from 2030 to 2050 is likely to involve substantial changes in the location and magnitude of generation and demand, increasing the need for timely transmission reinforcement and flexibility investment. Robust long term transmission planning is therefore important for assessing future network requirements and supporting strategic investment decisions.

The GB context provides a useful setting for examining these challenges. However, the underlying planning issues are not unique to GB. Many electricity systems are facing similar challenges arising from increasing renewable generation, electrification, uncertainty in future energy pathways and long lead times for network investment. The GB system is therefore used in this paper as a case study through which a more general planning framework can be developed and demonstrated.

For long term energy system planning, the National Energy System Operator (NESO), publishes the annual Future Energy Scenarios (FES)~\cite{NESO_FES2025}. FES provides four alternative pathways to 2050, reflecting different assumptions about technology deployment, consumer behaviour and policy development in Great Britain. By modelling interactions across electricity, heating, transport, hydrogen and other energy vectors, FES provides a consistent basis for exploring how the GB energy system could evolve under different future conditions.

However, FES does not specify how the transmission network should develop to accommodate these changes. The Electricity Ten Year Statement (ETYS)~\cite{ETYS} provides information on future transmission capability requirements across major GB network boundaries, but does not provide a detailed network expansion plan or specify the timing and location of individual reinforcements. Consequently, neither source provides a network model that can be directly used to study the long term transmission investment decisions required under alternative energy pathways. There is therefore a need to connect national energy system scenarios with network-level investment planning.

This paper addresses this gap by developing a multi-year stochastic framework for transmission expansion and energy storage planning under alternative future energy pathways. The framework links national energy transition scenarios to network-constrained investment decisions using an enhanced representation of the GB transmission system. Annual demand and generation trajectories are derived from FES. Transmission and energy storage investments are then optimised jointly over an annual planning horizon from 2030 to 2050 using a mixed-integer stochastic DC optimal power flow formulation.

The main contributions of the paper are:

\begin{itemize}
\item \textbf{Multi-year stochastic transmission and storage planning model.} A planning framework is developed to determine the scale, timing and location of transmission and storage investments under uncertainty, while explicitly representing investment lead times and commissioning periods.
\item \textbf{Linking energy pathways to network investment}. Industry data sources, including FES and ETYS, are integrated with an enhanced representative GB transmission network to translate national energy transition pathways into network-level investment requirements.
\item \textbf{Quantifying the value of coordinated transmission and storage investment}. In the GB electricity market, transmission and energy storage investment decisions are made by different actors and are not centrally co-optimised. Transmission investment requires forward-looking network planning, while storage investment is undertaken by market participants responding to market conditions and investment incentives. The joint optimisation developed in this paper therefore does not represent a practical centralised investment decision. Instead, it provides a counterfactual planning benchmark that quantifies the potential system value of better coordination between transmission development and storage deployment. This benchmark can provide evidence to inform market and policy arrangements that improve investment signals and enable the market to deliver storage where and when it provides the greatest system value.
\item \textbf{Assessment of alternative futures and planning strategies}. Different energy pathways and planning strategies are evaluated in terms of investment requirements, investment timing and location, system cost and supply adequacy, providing insight into how uncertainty in the energy transition affects long term transmission and storage planning.

\end{itemize}


The planning framework developed in this paper is intended for strategic long term planning rather than project-level transmission design. By bringing together pathway uncertainty, network constraints, transmission and storage investment, and asset-specific investment lead times and commissioning periods, it provides a consistent basis for examining how alternative energy futures translate into the scale, location and timing of transmission and flexibility investment.

Although the analysis focuses on the GB electricity system, the framework addresses challenges that are relevant to other power systems undergoing large-scale renewable integration and electrification. In particular, the joint planning formulation provides a useful benchmark for understanding the potential value of better alignment between network investment and the development of flexible resources such as energy storage, even where these decisions are made by different market participants and are not currently co-optimised.
\section{Literature Review}

\subsection{Future Energy Scenario (FES)}
The National Energy System Operator (NESO) publishes its Future Energy Scenarios (FES) each year, providing an outlook for how the energy system in Great Britain (GB) could develop towards 2050 and beyond~\cite{cao2025stochastic,guillen2025long}. FES considers a range of factors that may shape the future energy system, including electricity demand, renewable generation, electrification, hydrogen use and consumer behaviour. Based on different assumptions about these factors, FES 2025 presents four pathways: Holistic Transition (HT), Electric Engagement (EE), Hydrogen Evolution (HE), and Falling Behind (FB). Their main characteristics are summarized in Table~\ref{tab:FES_scenarios}.

\begin{table}
	\centering
	\caption{Summary of the four FES 2025 pathways~\cite{NESO_FES2025}.}
	\label{tab:FES_scenarios}
	\begin{tabular}{p{0.15\columnwidth} p{0.82\columnwidth}}
	\toprule
		\textbf{Pathway} & \textbf{Description} \\
	\midrule
		HT &
		A coordinated pathway towards net zero involving changes across multiple energy sectors, with a broad mix of electrification, hydrogen, consumer engagement, and energy-efficiency measures. \\
	\midrule		
EE &
		A highly electrified pathway in which consumers actively adopt low-carbon technologies, with extensive electrification of transport and heating and increasing reliance on renewable electricity. \\
			\midrule
HE &
		A pathway characterized by a greater role for hydrogen in decarbonizing energy demand, particularly in sectors where direct electrification is more difficult. \\
			\midrule
FB &
		A slower-transition pathway characterized by weaker progress in decarbonisation and lower adoption of low-carbon technologies, resulting in the net-zero target not being achieved by 2050. \\
		
\bottomrule
	\end{tabular}
\end{table}

Previous studies have widely employed FES to investigate the implications of the energy transition in GB. These studies include modeling future heating and electricity demand~\cite{watson2019decarbonising,watson2023predicting,simpson2026air}, assessing renewable energy integration and system flexibility~\cite{daymodelling,gill2024far,chyong2025zonal,zannane2025climate}, and examining the effects of transport electrification and sector coupling~\cite{dik2025sustainable,wang2025probabilistic,gillich2026end,zemani2025natural,kutlu2025comparative}. By providing a consistent set of assumptions across different energy sectors and planning horizons, FES provides an established framework for analyzing long term decarbonisation pathways in GB~\cite{brighty2025impacts,hutchinson2025qualitative,chen2025revealing}. Its widespread use also supports its adoption as a basis for long term power system analysis.

However, existing FES based studies are mainly suited to analyzing energy transitions at the national level. They provide limited guidance on where and when transmission network reinforcement should be implemented. FES describes how the energy system may develop under different future pathways, but does not directly determine the transmission investments required to accommodate these changes.

Therefore, although FES provides a comprehensive representation of possible future energy pathways, additional transmission network modeling is required before these scenarios can be used for long term Transmission Expansion Planning (TEP). Such modeling is needed to translate the changes in generation, demand and energy flows described by FES into specific transmission reinforcement decisions.

\subsection{Transmission Network Models for Great Britain}

To investigate the operation and future development of the GB transmission system, various network models have been developed with different levels of spatial and electrical detail. These models support power flow analysis, market simulation, long-term planning \cite{chaudry2023simulating,shafiekhani2025addressing}.
Detailed transmission network models represent major physical components, including substations, transmission lines, and transformers. Their ability to capture realistic power flows and network constraints makes them suitable for operational and security studies \cite{cao2025review,ding2016adjustable}. However, their scale and computational complexity can limit their application to long-term planning studies involving numerous scenarios and planning horizons.

Reduced transmission network models improve computational efficiency by simplifying network topology \cite{NESOREDUCEDMODEL}. They therefore provide a practical balance between network fidelity and computational tractability and are widely used in electricity market, transmission planning, and energy system optimisation studies \cite{gerber20122020,rabbani2014modelling,cao2025net}.
GB36 is a widely adopted reduced model of the GB transmission system due to its simple topology and representation of major transmission boundaries \cite{cao2025stochastic,azizipanah2024integration,mustapha2024novel}. However, its original configuration does not reflect recent network developments, including the Western HVDC Link and updated transmission capabilities reported in the ETYS \cite{ETYS}. Its implementation in the commercial DIgSILENT PowerFactory environment also limits accessibility and compatibility with open-source platforms such as MATPOWER, constraining its use in reproducible optimisation-based studies \cite{shafiekhani2025addressing,lyden2024pypsa}.

Although GB36 has been widely and successfully applied in research on GB power systems, there is still a need to update its line parameters and incorporate new lines added over the past decade or so. Enhancing the model’s compatibility with open-source optimisation frameworks is also crucial for supporting reproducible research and developing customized planning methods. This will enable it to address future challenges in the GB power system within the context of FES.

\subsection{Power Flow Modelling in Investment Studies}

Investment studies, in particular TEP, involves binary decisions that may model the timing, capacity and location of investment in power systems. Being combinatorial optimisation problems, these are often result in very large scale mixed integer linear programming problems~\cite{7769214,9087776,https://doi.org/10.1049/rpg2.12762}. The complexity of the problem increases when uncertainty and non-linearity of power flow is introduced. 


To maintain computational tractability, most long-term transmission and energy storage planning models employ DCOPF to represent network operation. Compared with ACOPF, its linear formulation is more suitable for large-scale problems involving multiple planning years, operating scenarios, and binary investment decisions. It also facilitates the formulation of MILP models and the application of decomposition techniques such as Benders decomposition \cite{cao2025review,9395702,6859489,https://doi.org/10.1049/joe.2018.9350}.

However, due to computational constraints and the focus of the models on specific subjects, there remain several significant limitations in the existing literature. Many studies simplify the planning horizon by selecting only a few representative planning years (such as 2030, 2040 and 2050) or a few ‘snapshots’ spanning several years \cite{shafiekhani2025addressing,lyden2024pypsa,watson2019decarbonising}. Although this approach reduces computational effort, it may produce discontinuous investment decisions and fail to capture gradual system evolution. Moreover, transmission and energy storage investments differ substantially in their development lead times and operational lifetimes: transmission projects typically require longer planning, approval, and construction periods, whereas storage can be deployed more rapidly. However, these differences are often neglected by assuming identical investment schedules \cite{lyden2024pypsa,shafiekhani2025addressing}. The economic representation adopted in many planning models remains relatively simplified, with practical investment factors such as inflation, financing costs, and time-dependent capital expenditure often represented through fixed investment costs or constant discount rates \cite{shafiekhani2025addressing,9395702}. Finally, the majority of planning methodologies continue to be validated using synthetic benchmark systems, such as the IEEE 24-bus, IEEE 118-bus, and IEEE 300-bus networks, which provide valuable platforms for methodological development but cannot adequately represent the geographical characteristics, network topology, and regional transmission constraints of real-world power systems \cite{10015767,7944676,9238453,7769214,9087776}. Consequently, the applicability of these planning approaches to realistic national transmission systems, such as the GB network under future FES pathways, remains insufficiently explored.

\subsection{Research Gap and Motivation}

Despite substantial progress in future energy scenario analysis, network modeling, and transmission and energy storage planning, these areas remain largely disconnected. In particular, few studies translate FES pathways and their regional demand and generation changes into long-term transmission investment decisions.
Reduced network models such as GB36 offer a practical balance between network detail and computational efficiency. However, existing implementations do not fully reflect recent infrastructure developments and transmission capability assumptions, while their reliance on proprietary software can limit interoperability and reproducibility.
Long-term planning studies also commonly rely on a limited number of planning years, which may produce discontinuous investment decisions and overlook differences in the lead times and operational lifetimes of transmission and storage assets. Moreover, many methods are validated only on synthetic test systems that cannot represent the geographical structure and regional transfer constraints of the GB network. Therefore, an integrated framework combining FES pathways, an updated GB network model, and realistic multi-year transmission and storage planning remains lacking.

\begin{table*}\small
	\centering
	\caption{Comparison of representative studies on future transmission planning for the Great Britain power system. T \& S in the last column represents studies that co-optimise planning of transmission and storage.}
	\label{tab:literature_review}
	\renewcommand{\arraystretch}{1.1}
	\begin{tabular}{lccccccc}
		\hline
		\rot{\textbf{Reference}} &
		\rot{\textbf{FES-based}} &
		\rot{\textbf{GB Network}} &
		\rot{\textbf{TEP}} &
		\rot{\textbf{Storage}} &
		\rot{\textbf{Annual Planning}} &
		\rot{\textbf{Stochastic}} &
		\rot{\textbf{T \& S }} \\
		\hline
		\cite{cao2025stochastic} &
		$\checkmark$ &
		$\checkmark$ &
		$\checkmark$ &
		$\times$ &
		$\times$ &
		$\checkmark$ &
		$\times$ \\
		\cite{watson2019decarbonising} &
		$\checkmark$ &
		$\checkmark$ &
		$\times$ &
		$\times$ &
		$\times$ &
		$\times$ &
		$\times$ \\
		\cite{simpson2026air} &
		$\checkmark$ &
		$\checkmark$ &
		$\times$ &
		$\times$ &
		$\checkmark$ &
		$\times$ &
		$\times$ \\
		\cite{dik2025sustainable} &
		$\checkmark$ &
		$\times$ &
		$\times$ &
		$\times$ &
		$\times$ &
		$\checkmark$ &
		$\checkmark$ \\
		\cite{wang2025probabilistic} &
		$\checkmark$ &
		$\times$ &
		$\times$ &
		$\checkmark$ &
		$\checkmark$ &
		$\checkmark$ &
		$\checkmark$ \\
		\cite{zemani2025natural} &
		$\checkmark$ &
		$\times$ &
		$\times$ &
		$\checkmark$ &
		$\times$ &
		$\times$ &
		$\checkmark$ \\
		\cite{kutlu2025comparative} &
		$\checkmark$ &
		$\times$ &
		$\times$ &
		$\checkmark$ &
		$\checkmark$ &
		$\checkmark$ &
		$\checkmark$ \\
		\cite{shafiekhani2025addressing} &
		$\checkmark$ &
		$\checkmark$ &
		$\checkmark$ &
		$\checkmark$ &
		$\times$ &
		$\times$ &
		$\checkmark$ \\
		\cite{ding2016adjustable} &
		$\times$ &
		$\times$ &
		$\times$ &
		$\times$ &
		$\checkmark$ &
		$\checkmark$ &
		$\times$ \\
		\cite{gerber20122020} &
		$\times$ &
		$\checkmark$ &
		$\times$ &
		$\times$ &
		$\checkmark$ &
		$\times$ &
		$\times$ \\
		\cite{rabbani2014modelling} &
		$\checkmark$ &
		$\checkmark$ &
		$\times$ &
		$\times$ &
		$\checkmark$ &
		$\times$ &
		$\times$ \\
		\cite{cao2025net} &
		$\checkmark$ &
		$\checkmark$ &
		$\checkmark$ &
		$\times$ &
		$\checkmark$ &
		$\checkmark$ &
		$\times$ \\
		\cite{azizipanah2024integration} &
		$\times$ &
		$\checkmark$ &
		$\times$ &
		$\checkmark$ &
		$\checkmark$ &
		$\times$ &
		$\checkmark$ \\
		\cite{mustapha2024novel} &
		$\checkmark$ &
		$\checkmark$ &
		$\times$ &
		$\times$ &
		$\checkmark$ &
		$\times$ &
		$\checkmark$ \\
		\cite{lyden2024pypsa} &
		$\checkmark$ &
		$\checkmark$ &
		$\times$ &
		$\times$ &
		$\times$ &
		$\checkmark$ &
		$\times$ \\
		\cite{7769214} &
		$\times$ &
		$\times$ &
		$\checkmark$ &
		$\times$ &
		$\checkmark$ &
		$\checkmark$ &
		$\checkmark$ \\
		\cite{9087776} &
		$\times$ &
		$\times$ &
		$\checkmark$ &
		$\times$ &
		$\checkmark$ &
		$\checkmark$ &
		$\times$ \\
		\cite{https://doi.org/10.1049/rpg2.12762} &
		$\times$ &
		$\times$ &
		$\checkmark$ &
		$\times$ &
		$\checkmark$ &
		$\checkmark$ &
		$\checkmark$ \\
		\cite{https://doi.org/10.1049/joe.2018.9350} &
		$\times$ &
		$\times$ &
		$\times$ &
		$\checkmark$ &
		$\checkmark$ &
		$\checkmark$ &
		$\checkmark$ \\
		\cite{6859489} &
		$\times$ &
		$\times$ &
		$\times$ &
		$\checkmark$ &
		$\checkmark$ &
		$\checkmark$ &
		$\checkmark$ \\
		\cite{qiu2016stochastic} &
		$\times$ &
		$\times$ &
		$\checkmark$ &
		$\checkmark$ &
		$\checkmark$ &
		$\checkmark$ &
		$\checkmark$ \\
		\cite{10015767} &
		$\times$ &
		$\times$ &
		$\checkmark$ &
		$\checkmark$ &
		$\checkmark$ &
		$\checkmark$ &
		$\checkmark$ \\
		\cite{7944676} &
		$\times$ &
		$\times$ &
		$\checkmark$ &
		$\times$ &
		$\checkmark$ &
		$\checkmark$ &
		$\checkmark$ \\
		\cite{9238453} &
		$\times$ &
		$\times$ &
		$\checkmark$ &
		$\times$ &
		$\checkmark$ &
		$\checkmark$ &
		$\checkmark$ \\
		\cite{9395702} &
		$\times$ &
		$\times$ &
		$\checkmark$ &
		$\times$ &
		$\checkmark$ &
		$\checkmark$ &
		$\checkmark$ \\
		\cite{8253865} &
		$\times$ &
		$\times$ &
		$\checkmark$ &
		$\times$ &
		$\checkmark$ &
		$\checkmark$ &
		$\checkmark$ \\
		\textbf{This paper} &
		$\checkmark$ &
		$\checkmark$ &
		$\checkmark$ &
		$\checkmark$ &
		$\checkmark$ &
		$\checkmark$ &
		$\checkmark$ \\
		\hline
	\end{tabular}
\end{table*}

Table~\ref{tab:literature_review} compares representative studies on future GB transmission planning and highlights the lack of an integrated framework combining updated GB network representation, FES pathways, and annual multi-year planning.
To address this gap, this paper integrates the FES 2025 pathways with an enhanced GB36 model for annual transmission and energy storage co-planning from 2030 to 2050 using stochastic DCOPF. It also distinguishes their commissioning times by incorporating line-specific preparation periods for transmission reinforcement while allowing storage to operate from the investment year.

\section{Methodology}
\subsection{Enhanced Great Britain 36-Bus Power System Model}

The original GB36 model developed by the system operator provides a reduced representation of the Great Britain transmission system in which major geographical regions and transmission corridors are represented by 36 buses and their associated interregional connections. However, the original model reflects an earlier configuration of the GB transmission network and does not include several major HVDC transmission developments that are relevant to the network configuration at the beginning of the planning horizon considered in this study. The original GB36 model is also based on interconnector development assumptions made in 2012, some of which have subsequently been cancelled or revised.

\begin{figure}[!t]
	\centering
	
	\subfloat[Original GB36 model \cite{NESO_GB_Reduced_Model}.]{
		\includegraphics[width=0.475\columnwidth]{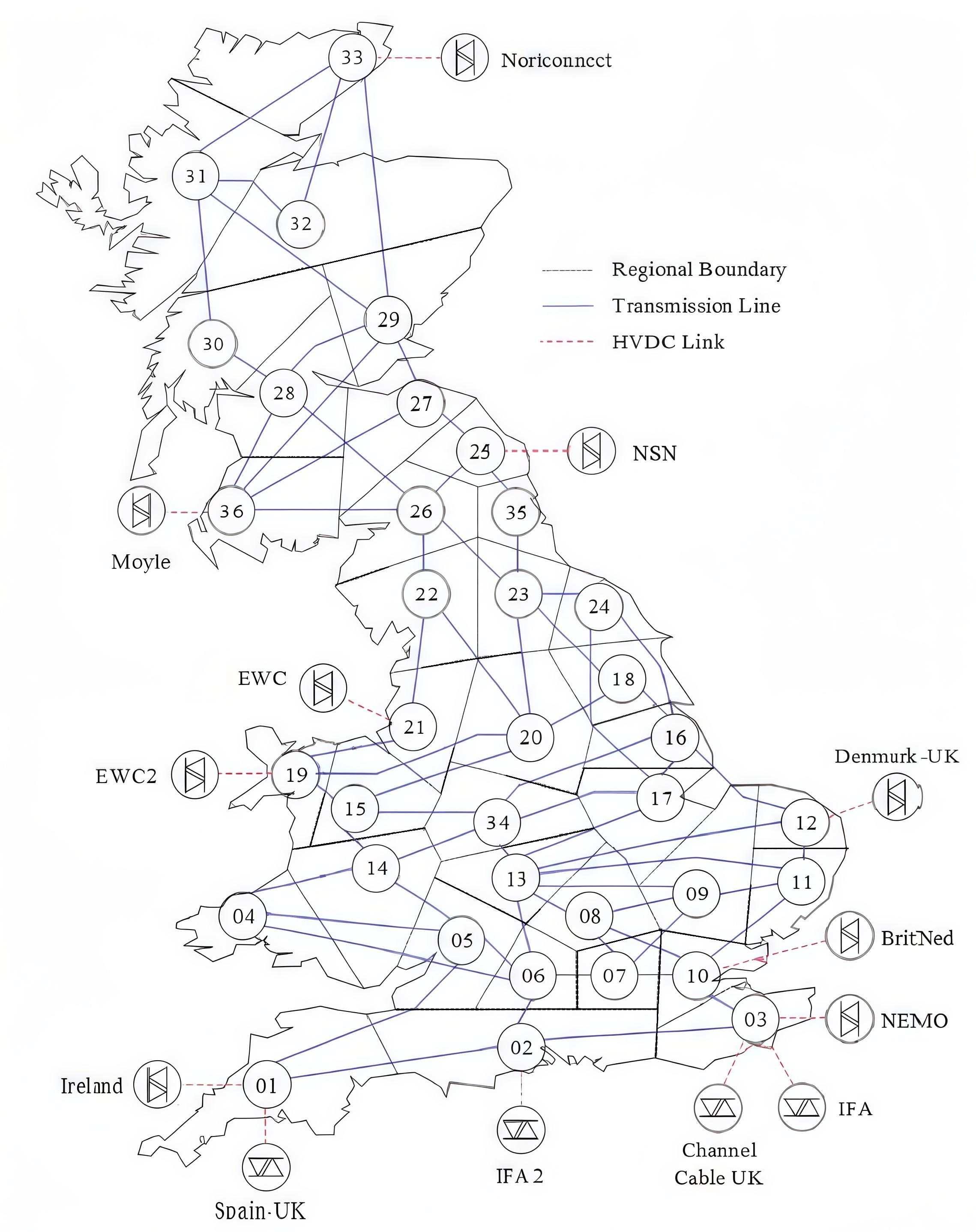}
		\label{fig:GB36_original}
	}
	\hfill
	\subfloat[Enhanced GB36 model.]{
		\includegraphics[width=0.475\columnwidth]{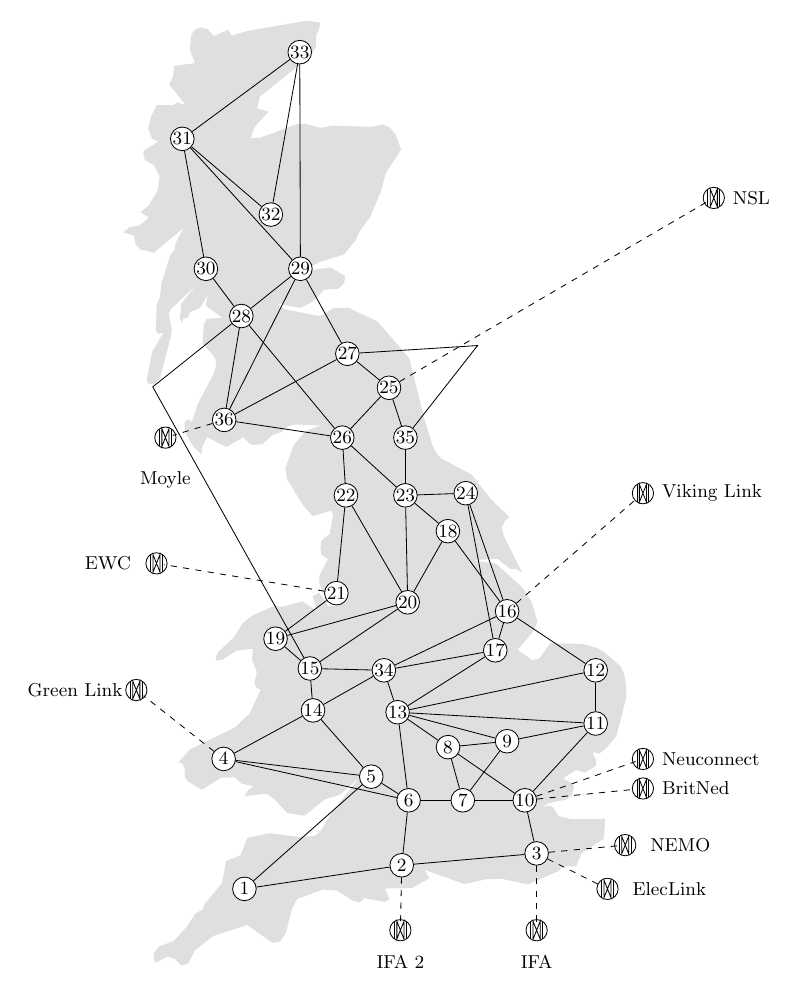}
		\label{fig:GB36_enhanced}
	}
	
	\caption{Comparison of the original and enhanced GB36 transmission network models:
		(a) original GB36 model from the electricity system operator, NESO, based on 2012 data; and (b) enhanced GB36 model incorporating recent HVDC links, and revised inter area capacities, representing a 2026 test case.}
	\label{fig:GB36_comparison}
\end{figure}

To provide a more up-to-date representation of the GB transmission system, the original GB36 topology is enhanced by incorporating two major HVDC links. The first is the Western HVDC Link, represented as an additional connection between buses 28 and 15, which has been in operation since 2019. The second is Eastern Green Link 1 (EGL1), represented between buses 27 and 35, which is scheduled to become operational in 2029. Since the planning horizon of this study begins in 2030, both links are therefore treated as existing transmission infrastructure in the initial network configuration. The enhanced GB36 model retains 36 buses while increasing the number of transmission corridors from 66 to 68. Other proposed EGL projects, including EGL2--EGL5, are not incorporated into the improved GB36 network. The resulting enhanced GB36 model is illustrated in Fig.~\ref{fig:GB36_comparison}(b).

The transmission capacities of the existing corridors are recalibrated using the latest boundary capability information reported in ETYS. Each GB36 transmission branch is assigned to its corresponding ETYS boundary according to its geographical location. Where multiple branches cross the same boundary, the published boundary transfer capability is assumed to be equally distributed among the associated branches. The recalibrated branch capacity is therefore calculated as

\begin{equation}
	F_{l}^{\max}
	=
	\frac{C_{b}^{\mathrm{ETYS}}}{N_{b}}, \forall l\in\mathcal{L}_{b}.
\end{equation}

\noindent
where $C_{b}^{\mathrm{ETYS}}$ denotes the transfer capability of transmission boundary $b$ reported in ETYS, $\mathcal{L}_{b}$ is the set of GB36 branches crossing boundary $b$, and $N_{b}=|\mathcal{L}_{b}|$ is the number of associated branches.

The generation and demand data used in this study are derived from the four scenarios presented in the FES 2025 report. For each scenario, the electricity demand is represented using the corresponding AES Peak values reported in FES 2025. In each planning year, system operation is represented by a representative 24-hour operating profile, enabling the model to capture the chronological variation in electricity demand and renewable generation while maintaining computational tractability.

During the planning period, both generating capacity and electricity demand are assumed to grow year by year, with the increments allocated according to each node’s initial share to maintain a constant spatial distribution. Since FES 2025 does not provide information on the locations or connection details of new interconnection lines, the interconnection capacity remains at the values set in FES 2025 throughout the planning period.

Each existing transmission line may be expanded to a maximum of twice its initial capacity. The upper limit for new energy storage capacity is set in accordance with the \textit{Electric Engagement} (EE) scenario of FES 2025, totalling 43.2 GW, and is distributed evenly across the 36 buses of GB36, i.e. a maximum of 1.2 GW per bus.

\subsection{Financial Assumptions and Investment Cost Modeling}

Load shedding is penalized at \pounds10,000/MWh to represent the high
economic cost of unserved energy and to discourage load shedding unless
required to maintain system feasibility.

All capital costs are expressed in 2030 constant prices and evaluated on a present-value basis. As transmission and storage assets may remain operational beyond the 2050 planning horizon, residual values are assigned to their remaining lifetimes.

Let $C_{a,y_0}^{\mathrm{CAPEX}}$ denote the reference capital cost of candidate
asset $a$ in the base year $y_0=2030$. Transmission capital costs are assumed
to remain constant in real terms, whereas storage capital costs decline at an
annual real rate $\gamma^{\mathrm{ES}}$ due to technological development and
learning effects. The year-specific undiscounted capital cost is therefore
defined as
\begin{equation}
	C_{a,y}^{\mathrm{CAPEX}}
	=
	\begin{cases}
		C_{l,y_0}^{\mathrm{TL}},
		& a=l,\quad l\in\mathcal{L},\\[1mm]
		C_{k,y_0}^{\mathrm{ES}}
		\left(1-\gamma^{\mathrm{ES}}\right)^{y-y_0},
		& a=k,\quad k\in\mathcal{K}.
	\end{cases}
\end{equation}

\noindent
where $\mathcal{L}$ and $\mathcal{K}$ denote the sets of transmission
lines and candidate energy storage sites, respectively.

The storage cost reduction rate represents technology-driven changes in real
capital cost and is therefore applied separately from the financial discount
rate. Accordingly, the present value of the capital expenditure associated
with an investment decision in year $y$ is

\begin{equation}
	PV_{a,y}^{\mathrm{CAPEX}}
	=
	\frac{
		C_{a,y}^{\mathrm{CAPEX}}
	}{
		(1+r)^{y-y_0}
	},
\end{equation}
\noindent
where $r$ denotes the real discount rate. A real discount rate of $4.5\%$ is
adopted consistently for transmission and energy storage investments.

A line specific preparation period represents the time required to plan, obtain approvals for, and construct each transmission reinforcement, measured from the investment decision to the commissioning of that reinforcement. Let $\tau_l^{\mathrm{prep}}$ denote the preparation period of candidate
transmission line $l$. We determine $\tau_l^{\mathrm{prep}}$ based on the length of the candidate line, using the following assumptions:

\begin{equation}
	\tau_l^{\mathrm{prep}}
	=
	\begin{cases}
		1, & d_l < 100\ \mathrm{km},\\
		2, & 100 \le d_l < 200\ \mathrm{km},\\
		3, & d_l \ge 200\ \mathrm{km},
	\end{cases}
\end{equation}
\noindent
where $d_l$ denotes the length of candidate transmission line $l$. The
resulting preparation periods are predetermined during the data preparation
stage and imported from the input dataset as fixed parameters.

Energy storage systems (ESS) are assumed to have no explicit preparation period and
therefore become operational in the year in which the investment is made. Based on these assumptions, for a general candidate asset $a$, the commissioning year is defined as:

\begin{equation}
	y_{a,y}^{\mathrm{com}}
	=
	\begin{cases}
		y+\tau_l^{\mathrm{prep}},
		& a=l,\quad l\in\mathcal{L},\\
		y,
		& a=k,\quad k\in\mathcal{K}.
	\end{cases}
\end{equation}

\noindent
where $y$ is the investment year and $\tau_l^{\mathrm{prep}}$ is the preparation period of candidate transmission line $l$.

Given the final planning year $Y=2050$, the number of operational years
represented within the planning horizon is

\begin{equation}
	n_{a,y}^{\mathrm{used}}
	=
	\max
	\left[
	0,\,
	Y-y_{a,y}^{\mathrm{com}}+1
	\right].
\end{equation}

The remaining lifetime of asset $a$ at the end of the planning horizon is

\begin{equation}
	n_{a,y}^{\mathrm{rem}}
	=
	\max
	\left[
	0,\,
	n_a-n_{a,y}^{\mathrm{used}}
	\right],
\end{equation}

\noindent
where $n_a$ denotes the assumed technical or economic lifetime of the
corresponding asset.

Technology-specific asset lifetimes are adopted in this study:

\begin{equation}
	n_a=
	\begin{cases}
		45, & a \in \mathcal{L}^{\mathrm{Transmission}},\\
		21, & a \in \mathcal{K}^{\mathrm{ES}},
	\end{cases}
\end{equation}
\noindent
where $\mathcal{L}^{\mathrm{Transmission}}$ denotes transmission assets and
$\mathcal{K}^{\mathrm{ES}}$ denotes energy storage assets.

Assuming straight-line depreciation over the lifetime of each asset, the
residual value of asset $a$ at the end of the planning horizon is calculated
using the capital cost applicable to the corresponding investment year:

\begin{equation}
	RV_{a,y}^{2050}
	=
	C_{a,y}^{\mathrm{CAPEX}}
	\frac{
		n_{a,y}^{\mathrm{rem}}
	}{
		n_a
	}.
\end{equation}

The residual value is subsequently discounted from 2050 to the 2030 base year:

\begin{equation}
	PV_{a,y}^{RV}
	=
	\frac{
		RV_{a,y}^{2050}
	}{
		(1+r)^{Y-y_0}
	}.
\end{equation}

Consequently, the effective present-value investment cost associated with
constructing asset $a$ in year $y$ is defined as

\begin{equation}
	C_{a,y}^{\mathrm{eff}}
	=
	PV_{a,y}^{\mathrm{CAPEX}}
	-
	PV_{a,y}^{RV}.
\end{equation}

Substituting the previous expressions gives

\begin{equation}
	C_{a,y}^{\mathrm{eff}}
	=
	\frac{
		C_{a,y}^{\mathrm{CAPEX}}
	}{
		(1+r)^{y-y_0}
	}
	-
	\frac{
		C_{a,y}^{\mathrm{CAPEX}}
	}{
		(1+r)^{Y-y_0}
	}
	\left(
	\frac{
		n_{a,y}^{\mathrm{rem}}
	}{
		n_a
	}
	\right).
\end{equation}

Effective investment costs are pre calculated for each candidate asset and each possible investment year before solving the transmission expansion model. These costs are then passed to the model as fixed parameters, so that the optimisation does not need to explicitly model the individual components of the investment cost calculation.

The pre calculated effective investment cost accounts for the timing of the investment, asset specific preparation periods and lifetimes, and the resulting period over which the asset is operational within the planning horizon. For energy storage systems, the calculation also accounts for assumed future reductions in storage costs. The costs are further adjusted for the relevant discounting and residual value of the asset at the end of the planning horizon. Thus, the model uses a single effective investment cost for each candidate asset and investment year that captures the relevant timing and lifetime effects over the planning horizon.

\subsection{Decision Variables}

\subsubsection{Continuous Variables}

\begin{align*}
    &P^{G}_{g,y,h,s} \ge 0
        &&\text{Conventional generator output};\\
    &P^{W}_{w,y,h,s} \ge 0
        &&\text{Utilized wind power};\\
    &P^{PV}_{p,y,h,s} \ge 0
        &&\text{Utilized PV power};\\
    &P^{L}_{l,y,h,s} \in \mathbb{R}
        &&\text{Power flow on line } l;\\
    &\theta_{b,y,h,s} \in \mathbb{R}
        &&\text{Voltage phase angle};\\
    &P^{S}_{k,y,h,s} \in \mathbb{R}
        &&\text{Net storage injection, positive when discharging};\\
    &P^{S,\mathrm{in}}_{k,y,h,s} \ge 0
        &&\text{Charging power};\\
    &P^{S,\mathrm{out}}_{k,y,h,s} \ge 0
        &&\text{Discharging power};\\
    &E_{k,y,h,s} \ge 0
        &&\text{State of charge of storage unit } k;\\
    &P^{\mathrm{in}}_{b,i,y,h,s} \ge 0
        &&\text{Imported power};\\
    &P^{\mathrm{out}}_{b,i,y,h,s} \ge 0
        &&\text{Exported power};\\
    &L^{\mathrm{shed}}_{b,y,h,s} \ge 0
        &&\text{Load shedding at bus } b.
\end{align*}

\subsubsection{Binary Variables}

\begin{align*}
&x^{\mathrm{build}}_{l,y}\in\{0,1\}
&&\text{1 if investment in line }l\text{ is initiated in year }y;\\
&x^{\mathrm{op}}_{l,y}\in\{0,1\}
&&\text{1 if line reinforcement is operational;}\\
&z^{\mathrm{build}}_{k,y}\in\{0,1\}
&&\text{1 if storage }k\text{ is built in year }y;\\
&z^{\mathrm{op}}_{k,y}\in\{0,1\}
&&\text{1 if storage }k\text{ is operational;}\\
&u^{\mathrm{in}}_{k,y,h,s}\in\{0,1\}
&&\text{1 if storage }k\text{ is charging;}\\
&u^{\mathrm{out}}_{k,y,h,s}\in\{0,1\}
&&\text{1 if storage }k\text{ is discharging.}
\end{align*}

TEP decisions are represented by binary variables for each candidate corridor and planning year. Reinforcement capacity becomes available only after the specified preparation period, allowing the model to optimize both investment selection and timing.
ESS investment variables determine the location and capacity of storage at candidate buses. Storage becomes operational in the investment year, with its hourly operation constrained by capacity limits, state-of-charge dynamics, efficiency, and charging and discharging limits.

\subsection{Objective Function}

The optimisation objective minimizes the total system cost over the planning
horizon, including generation cost, interconnector trading cost, renewable
curtailment cost, storage operating cost, load-shedding penalty, and the
effective present-value investment costs of transmission and energy storage
assets.

\begin{align}
\min \ \ &
365
\sum_{y\in\mathcal{Y}}
\sum_{h\in\mathcal{H}}
\sum_{s\in\Omega}
\Bigg[
\sum_{g\in\mathcal{G}}
\pi_s c^G_{g,y,s}P^G_{g,y,h,s}
+
\sum_{(b,i)\in IntBus}
\pi_s
\left(
c^{\mathrm{in}}_iP^{\mathrm{in}}_{b,i,y,h,s}
-
c^{\mathrm{out}}_iP^{\mathrm{out}}_{b,i,y,h,s}
\right)
\notag\\
&+
\sum_{w\in\mathcal{W}}
\pi_s c^{W}_{w,y,h,s}
\left(
\overline{P}^{W}_{w,y,h,s}
-
P^{W}_{w,y,h,s}
\right)
+
\sum_{p\in\mathcal{P}}
\pi_s c^{PV}_{p,y,h,s}
\left(
\overline{P}^{PV}_{p,y,h,s}
-
P^{PV}_{p,y,h,s}
\right)
\notag\\
&+
\sum_{k\in\mathcal{K}}
\pi_s
\left(
c^{\mathrm{ch}}_kP^{S,\mathrm{in}}_{k,y,h,s}
+
c^{\mathrm{dis}}_kP^{S,\mathrm{out}}_{k,y,h,s}
\right)
+
\sum_{b\in\mathcal{B}}
\pi_s C^{\mathrm{shed}}
L^{\mathrm{shed}}_{b,y,h,s}
\Bigg]
\notag\\
&+
\sum_{l\in\mathcal{L}}
\sum_{y\in\mathcal{Y}^{\mathrm{TL}}_l}
C_{l,y}^{\mathrm{TL,eff}}
x^{\mathrm{build}}_{l,y}
+
\sum_{k\in\mathcal{K}}
\sum_{y\in\mathcal{Y}}
C_{k,y}^{\mathrm{ES,eff}}
z^{\mathrm{build}}_{k,y}
\end{align}

\noindent
where $C_{l,y}^{\mathrm{TL,eff}}$ represents the effective present-value
investment cost of transmission line $l$ if the investment decision is made
in year $y$, while $C_{k,y}^{\mathrm{ES,eff}}$ represents the corresponding
effective present-value investment cost of energy storage asset $k$.

Because each transmission reinforcement has a line-specific preparation
period, its feasible investment-year set is defined as

\begin{equation}
{
\mathcal{Y}^{\mathrm{TL}}_l
=
\left\{
y\in\mathcal{Y}
\mid
y+\tau_l^{\mathrm{prep}}\le Y
\right\},
\qquad
\forall l\in\mathcal{L}.
}
\end{equation}

Therefore, a transmission investment is permitted only if the corresponding
reinforcement can become operational no later than the final planning year.
Energy storage may be constructed in any year throughout the planning horizon,
$\mathcal{Y}=\{2030,\ldots,2050\}$.

\subsection{Constraints}

\subsubsection{Power Balance}

For each $b,y,h,s$:

\begin{align}
&
\sum_{(b,g)\in GenBus} P^G_{g,y,h,s}
+
\sum_{(b,w)\in WinBus} P^{W}_{w,y,h,s}
+
\sum_{(b,p)\in SolBus} P^{PV}_{p,y,h,s}
\notag\\
&
+
\sum_{(b,k)\in StoBus}
P^{S,\mathrm{out}}_{k,y,h,s}
+
\sum_{\substack{(b',i)\in IntBus\\b'=b}}
P^{\mathrm{in}}_{b',i,y,h,s}
\notag\\
&
-
\sum_{\substack{l\in\mathcal{L}\\o(l)=b}}
P^{L}_{l,y,h,s}
+
\sum_{\substack{l\in\mathcal{L}\\r(l)=b}}
P^{L}_{l,y,h,s}
\notag\\
&=
\sum_{(b,d)\in DemBus}
D_{d,y,h,s}
-
L^{\mathrm{shed}}_{b,y,h,s}
\notag\\
&
+
\sum_{(b,k)\in StoBus}
P^{S,\mathrm{in}}_{k,y,h,s}
+
\sum_{\substack{(b',i)\in IntBus\\b'=b}}
P^{\mathrm{out}}_{b',i,y,h,s}.
\end{align}

\subsubsection{Reference Bus Angle}

\begin{equation}
\theta_{b^{\mathrm{ref}},y,h,s}=0,
\qquad
\forall y,h,s.
\end{equation}

\subsubsection{DC Power Flow and Angle Difference Constraints}

\begin{equation}
P^{L}_{l,y,h,s}
=
-B_l
\left(
\theta_{o(l),y,h,s}
-
\theta_{r(l),y,h,s}
\right),
\qquad
\forall l,y,h,s.
\end{equation}

\begin{equation}
\underline{\Delta\theta}
\le
\theta_{o(l),y,h,s}
-
\theta_{r(l),y,h,s}
\le
\overline{\Delta\theta},
\qquad
\forall l,y,h,s.
\end{equation}

\subsubsection{Line Capacity and Expansion Logic}

\begin{equation}
-
\left(
\overline{S}_l+\Delta S_lx^{\mathrm{op}}_{l,y}
\right)
\le
P^{L}_{l,y,h,s}
\le
\overline{S}_l+\Delta S_lx^{\mathrm{op}}_{l,y},
\qquad
\forall l,y,h,s.
\end{equation}

\subsubsection{Transmission Preparation and Operation Logic}

Transmission reinforcement is subject to a line-specific preparation period
$\tau_l^{\mathrm{prep}}$, representing the interval between the investment
decision and the year in which the reinforcement becomes operational.

The operational status of transmission reinforcement $l$ is linked directly
to all investment decisions whose preparation periods have been completed:

\begin{equation}
{
x^{\mathrm{op}}_{l,y}
=
\sum_{\substack{\tilde{y}\in\mathcal{Y}\\
\tilde{y}+\tau_l^{\mathrm{prep}}\le y}}
x^{\mathrm{build}}_{l,\tilde{y}},
\qquad
\forall l\in\mathcal{L},\ y\in\mathcal{Y}.
}
\end{equation}

Thus, an investment decision made in year $\tilde{y}$ becomes operational
only from year $\tilde{y}+\tau_l^{\mathrm{prep}}$ onward.

Each candidate transmission reinforcement may be initiated at most once:

\begin{equation}
{
\sum_{y\in\mathcal{Y}^{\mathrm{TL}}_l}
x^{\mathrm{build}}_{l,y}
\le
1,
\qquad
\forall l\in\mathcal{L}.
}
\end{equation}

Investment decisions whose preparation periods would extend beyond the
planning horizon are prohibited:

\begin{equation}
{
x^{\mathrm{build}}_{l,y}=0,
\qquad
\forall l\in\mathcal{L},\
y+\tau_l^{\mathrm{prep}}>Y.
}
\end{equation}

The preparation period is a predetermined line-specific parameter:

\begin{equation}
\tau_l^{\mathrm{prep}}\in\{1,2,3\},
\qquad
\forall l\in\mathcal{L}.
\end{equation}

\subsubsection{Storage Construction and Operation Logic}

Each candidate storage asset may be constructed at most once:

\begin{equation}
\sum_{y\in\mathcal{Y}}
z^{\mathrm{build}}_{k,y}
\le
1,
\qquad
\forall k\in\mathcal{K}.
\end{equation}

Since storage is assumed to become operational in the investment year, its
operational status is determined by

\begin{equation}
z^{\mathrm{op}}_{k,y}
=
\sum_{\tilde{y}\in\mathcal{Y}:\tilde{y}\le y}
z^{\mathrm{build}}_{k,\tilde{y}},
\qquad
\forall k\in\mathcal{K},\ y\in\mathcal{Y}.
\end{equation}

\subsubsection{Conventional Generation Limits and Ramping}

\begin{equation}
\underline{P}^{G}_{g,y,s}
\le
P^{G}_{g,y,h,s}
\le
\overline{P}^{G}_{g,y,s},
\qquad
\forall g,y,h,s.
\end{equation}

\begin{align}
P^{G}_{g,y,h+1,s}-P^{G}_{g,y,h,s}
&\le
RU_{g,y,s},
&&\forall g,y,s,\ h+1\in\mathcal{H},\\
P^{G}_{g,y,h+1,s}-P^{G}_{g,y,h,s}
&\ge
-RD_{g,y,s},
&&\forall g,y,s,\ h+1\in\mathcal{H}.
\end{align}

\subsubsection{Wind and PV Output Limits}

\begin{equation}
0
\le
P^{W}_{w,y,h,s}
\le
\overline{P}^{W}_{w,y,h,s},
\qquad
\forall w,y,h,s.
\end{equation}

\begin{equation}
0
\le
P^{PV}_{p,y,h,s}
\le
\overline{P}^{PV}_{p,y,h,s},
\qquad
\forall p,y,h,s.
\end{equation}

\subsubsection{Interconnector Capacity Constraints}

\begin{align}
0
\le
P^{\mathrm{in}}_{b,i,y,h,s}
&\le
\overline{P}^{\mathrm{in}}_i,
&&\forall (b,i)\in IntBus,\ y,h,s,\\
0
\le
P^{\mathrm{out}}_{b,i,y,h,s}
&\le
\overline{P}^{\mathrm{out}}_i,
&&\forall (b,i)\in IntBus,\ y,h,s.
\end{align}

\subsubsection{Storage Net Power Definition}

Storage net injection is positive when discharging:

\begin{equation}
P^{S}_{k,y,h,s}
=
P^{S,\mathrm{out}}_{k,y,h,s}
-
P^{S,\mathrm{in}}_{k,y,h,s},
\qquad
\forall k,y,h,s.
\end{equation}

\subsubsection{Storage Energy Balance}

For all $k\in\mathcal{K}$, $y\in\mathcal{Y}$, $s\in\Omega$, and
$h+1\in\mathcal{H}$:

\begin{equation}
E_{k,y,h+1,s}
=
E_{k,y,h,s}
+
\eta^{\mathrm{ch}}_k
P^{S,\mathrm{in}}_{k,y,h,s}
-
\frac{1}{\eta^{\mathrm{dis}}_k}
P^{S,\mathrm{out}}_{k,y,h,s}.
\end{equation}

The state of charge is bounded by the installed energy capacity:

\begin{equation}
0
\le
E_{k,y,h,s}
\le
\overline{E}_k z^{\mathrm{op}}_{k,y},
\qquad
\forall k,y,h,s.
\end{equation}

For representative-day operation, the storage energy level is assumed to be
cyclic:

\begin{equation}
E_{k,y,1,s}
=
E_{k,y,|\mathcal{H}|,s},
\qquad
\forall k,y,s.
\end{equation}

\subsubsection{Storage Power Limits and Mode Logic}

\begin{align}
P^{S,\mathrm{in}}_{k,y,h,s}
&\le
\overline{P}^{S}_k z^{\mathrm{op}}_{k,y},
&&\forall k,y,h,s,\\
P^{S,\mathrm{out}}_{k,y,h,s}
&\le
\overline{P}^{S}_k z^{\mathrm{op}}_{k,y},
&&\forall k,y,h,s,\\
-\overline{P}^{S}_k z^{\mathrm{op}}_{k,y}
\le
P^{S}_{k,y,h,s}
&\le
\overline{P}^{S}_k z^{\mathrm{op}}_{k,y},
&&\forall k,y,h,s.
\end{align}

The charging and discharging modes are mutually exclusive:

\begin{align}
u^{\mathrm{in}}_{k,y,h,s}
&\le
z^{\mathrm{op}}_{k,y},
&&\forall k,y,h,s,\\
u^{\mathrm{out}}_{k,y,h,s}
&\le
z^{\mathrm{op}}_{k,y},
&&\forall k,y,h,s,\\
u^{\mathrm{in}}_{k,y,h,s}
+
u^{\mathrm{out}}_{k,y,h,s}
&\le
z^{\mathrm{op}}_{k,y},
&&\forall k,y,h,s,\\
P^{S,\mathrm{in}}_{k,y,h,s}
&\le
\overline{P}^{S}_k u^{\mathrm{in}}_{k,y,h,s},
&&\forall k,y,h,s,\\
P^{S,\mathrm{out}}_{k,y,h,s}
&\le
\overline{P}^{S}_k u^{\mathrm{out}}_{k,y,h,s},
&&\forall k,y,h,s.
\end{align}

\subsubsection{Storage Ramping Constraints}

Let $M_k=\overline{P}^{S}_k$. For all $k,y,s$ and
$h+1\in\mathcal{H}$:

\begin{align}
P^{S,\mathrm{in}}_{k,y,h+1,s}
-
P^{S,\mathrm{in}}_{k,y,h,s}
&\le
RU^{S,\mathrm{ch}}_k
+
M_k(1-z^{\mathrm{op}}_{k,y}),
\\
P^{S,\mathrm{out}}_{k,y,h+1,s}
-
P^{S,\mathrm{out}}_{k,y,h,s}
&\le
RU^{S,\mathrm{dis}}_k
+
M_k(1-z^{\mathrm{op}}_{k,y}).
\end{align}

\section{Case Study}

The investment and operating cost assumptions are summarized in Table~\ref{tab:investment_cost}.
Capital costs are expressed in constant 2030 prices. Transmission investment
costs depend on additional transfer capacity and corridor length. The base-year
storage capital cost is \pounds300,000/MWh and declines over time as described
in Section~III. Investment costs are evaluated as effective present values,
accounting for discounting and the residual value at the end of the planning
horizon. The simulations use a real discount rate of 4.5\%, a transmission
asset lifetime of 45 years, and a storage asset lifetime of 21 years.
Transmission operating costs are assumed to be zero, while a cost of
\pounds30/MWh is applied separately to storage charging and discharging.

\begin{table}[H]
	\centering
	\caption{Investment and operating cost assumptions used in the simulations.}
	\label{tab:investment_cost}
	\begingroup
	\renewcommand{\arraystretch}{1.15}
	\setlength{\tabcolsep}{4pt}
	\resizebox{\columnwidth}{!}{%
	\begin{tabular}{lcc}
		\hline
		\textbf{Asset} &
		\textbf{CAPEX} &
		\textbf{OPEX} \\
		\hline
		
		AC overhead transmission line
		& \pounds570/MW/km
		& 0 \\
		
		Western HVDC Link
		& \pounds1,263/MW/km
		& 0 \\
		
		Eastern Green Link 1 (EGL1)
		& \pounds6,377/MW/km
		& 0 \\
		
		Energy storage
		& \pounds300,000/MWh
		& \pounds30/MWh (charge/discharge) \\
		\hline
	\end{tabular}%
	}
\endgroup
\end{table}
\subsection{Simulation Setup}

The optimisation problem is formulated as a mixed-integer linear programming (MILP) model based on stochastic DC optimal power flow and is solved using the Gurobi Optimizer. All simulations are performed on a workstation equipped with an Intel\textsuperscript{\textregistered} Core\texttrademark\ i7-9850H processor (2.60 GHz) and 32 GB RAM. The equivalent investment costs for transmission reinforcement and energy storage systems are calculated according to the financial assumptions presented in Section~III and are imported into the optimisation model as fixed input parameters.

Six simulation cases are considered to assess the value of coordinated transmission expansion planning (TEP) and energy storage system (ESS) investment under different assumptions about transmission preparation periods and future energy pathways, as summarised in Table~\ref{tab:case_setup}. The cases are designed to study the effects of (a) co optimising transmission and storage investment, (b) accounting for transmission preparation periods, and (c) considering multiple possible future energy pathways rather than planning against a single future scenario.

\begin{table}[H]
	\centering
	\caption{Summary of simulation cases.}
	\label{tab:case_setup}
	\renewcommand{\arraystretch}{1.2}
	\begin{tabular}{lcccc}
		\hline
		\textbf{Case} &
		\textbf{Transmission} &
		\textbf{Storage} &
		\textbf{Scenario} &
		\textbf{Preparation period} \\
		\hline
		
		A & $\checkmark$ & $\checkmark$ & Four FES & $\checkmark$\\
		
		B & $\checkmark$ & $\checkmark$ & Four FES & $\times$\\
		
		C & $\checkmark$ & $\times$ & Four FES & $\checkmark$\\
		
		D & $\times$ & $\checkmark$ & Four FES & $\checkmark$\\
		
		E & $\checkmark$ & $\checkmark$ & EE only & $\checkmark$\\
		
		F & $\checkmark$ & $\checkmark$ & FB only & $\checkmark$\\
		
		\hline
	\end{tabular}
\end{table}

Case A represents the proposed coordinated planning framework, in which TEP and ESS investments are co optimised within the stochastic planning model, while accounting for the line specific preparation period of each transmission reinforcement. This case provides the reference formulation against which the alternative investment strategies and planning assumptions are evaluated.

Case B uses the same coordinated TEP and ESS planning framework as Case A, but assumes that each transmission reinforcement becomes operational immediately in its investment year, thereby removing the transmission preparation period. Comparing Cases A and B isolates the effect of transmission preparation periods on investment decisions, asset deployment, and overall system costs.

Case C allows investment in transmission reinforcements but does not permit additional ESS investment. Conversely, Case D allows ESS expansion while keeping transmission capacities fixed throughout the planning horizon. Comparing Cases C and D with the coordinated Case A quantifies the value of jointly optimising transmission and storage investment relative to planning each asset class independently. These cases also indicate how the two technologies can substitute for or complement one another within the stochastic planning framework.

Cases E and F apply the coordinated TEP and ESS planning framework to the FES 2025 Electrification and Energy (EE) and Falling Behind (FB) scenarios, respectively. These cases examine how the optimal investment strategy changes under substantially different assumptions about the future energy system. In contrast to planning against a single future pathway with perfect foresight, considering multiple future pathways allows the investment strategy to account for uncertainty in the evolution of demand, generation, and other system conditions. The comparison therefore illustrates the implications of future pathway uncertainty for coordinated transmission and storage investment.

\subsection{Spatial Distribution of Expansion Decisions}
\label{subsubsec:spatial_expansion}

\begin{figure*}[!t]
	\centering
	
	\subfloat[Coordinated TEP and storage planning.]
	{
		\includegraphics[width=0.40\textwidth]
		{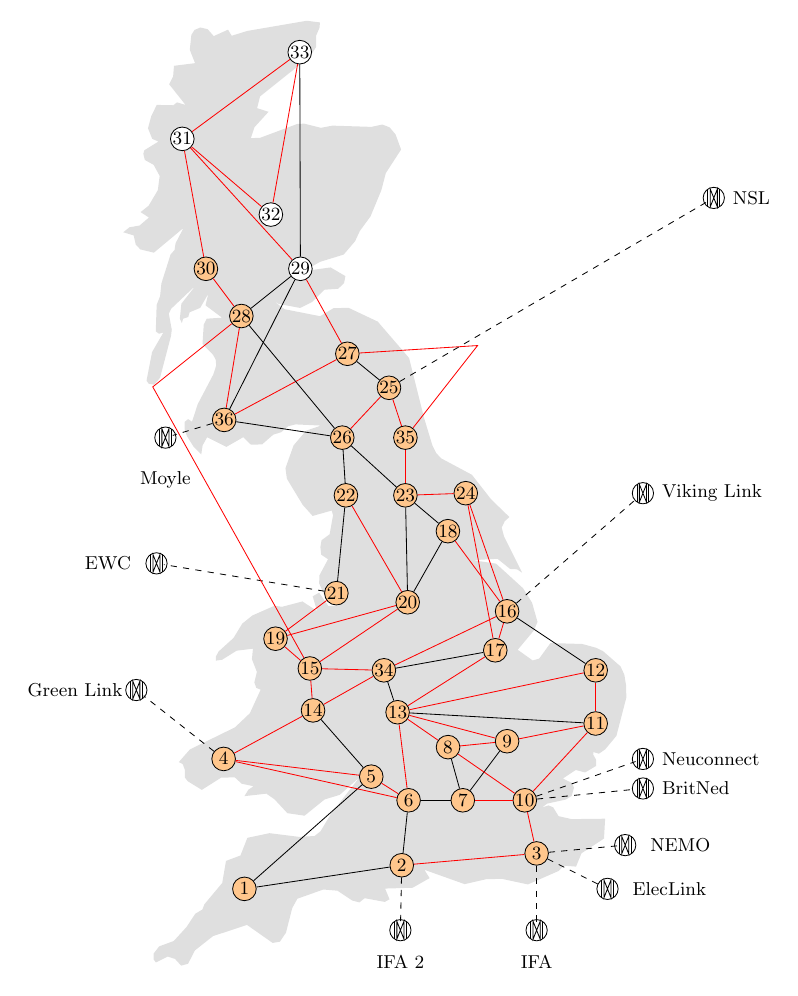}
		\label{fig:network_tep_storage}
	}
	\hfill
	\subfloat[Coordinated planning without time to build constraint.]
	{
		\includegraphics[width=0.40\textwidth]
		{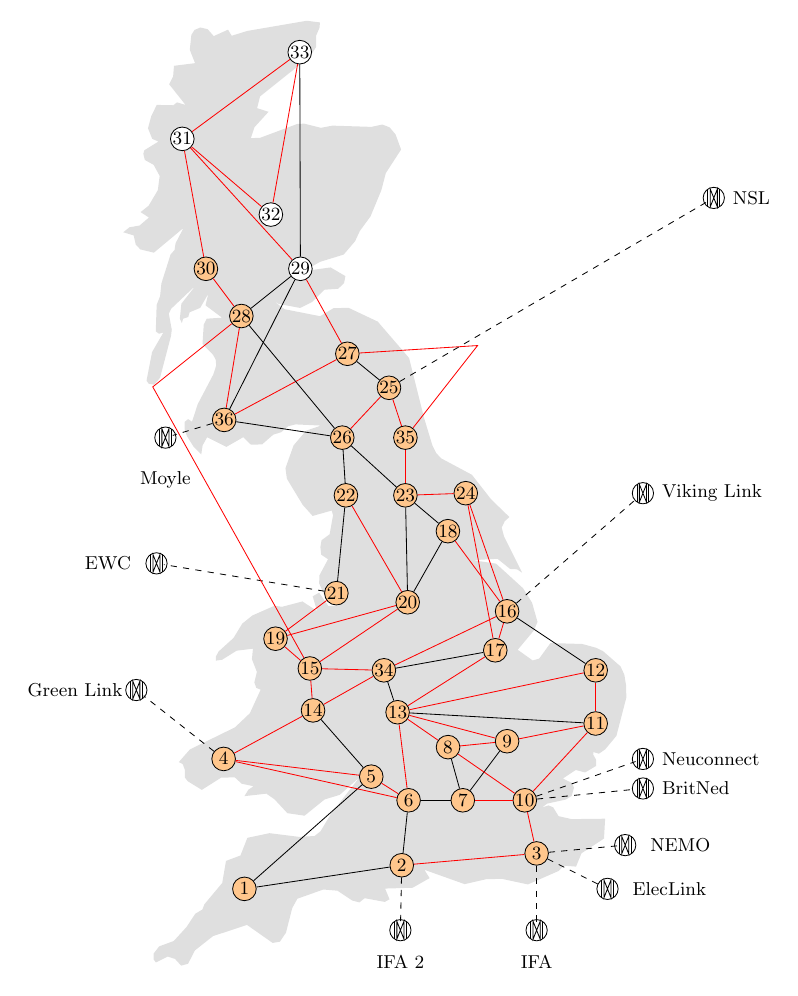}
		\label{fig:network_tep_storage_without_preparation}
	}
	
	\subfloat[Transmission-only planning.]
	{
		\includegraphics[width=0.40\textwidth]
		{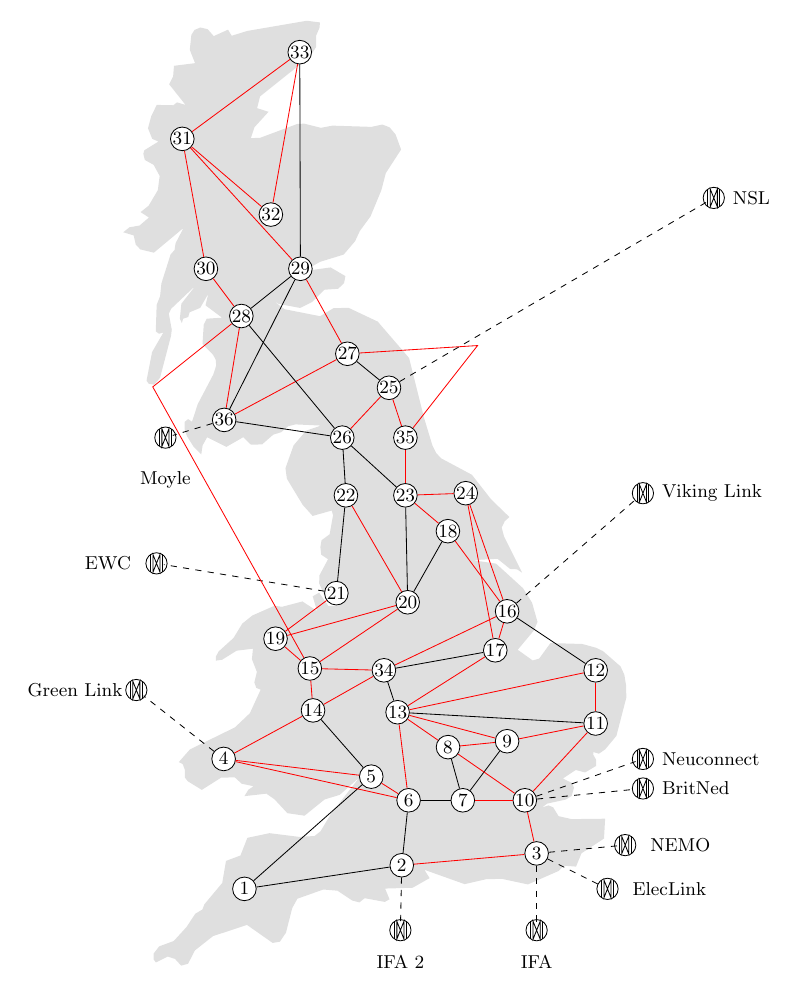}
		\label{fig:network_tep_only}
	}
	\hfill
	\subfloat[Storage-only planning.]
	{
		\includegraphics[width=0.40\textwidth]
		{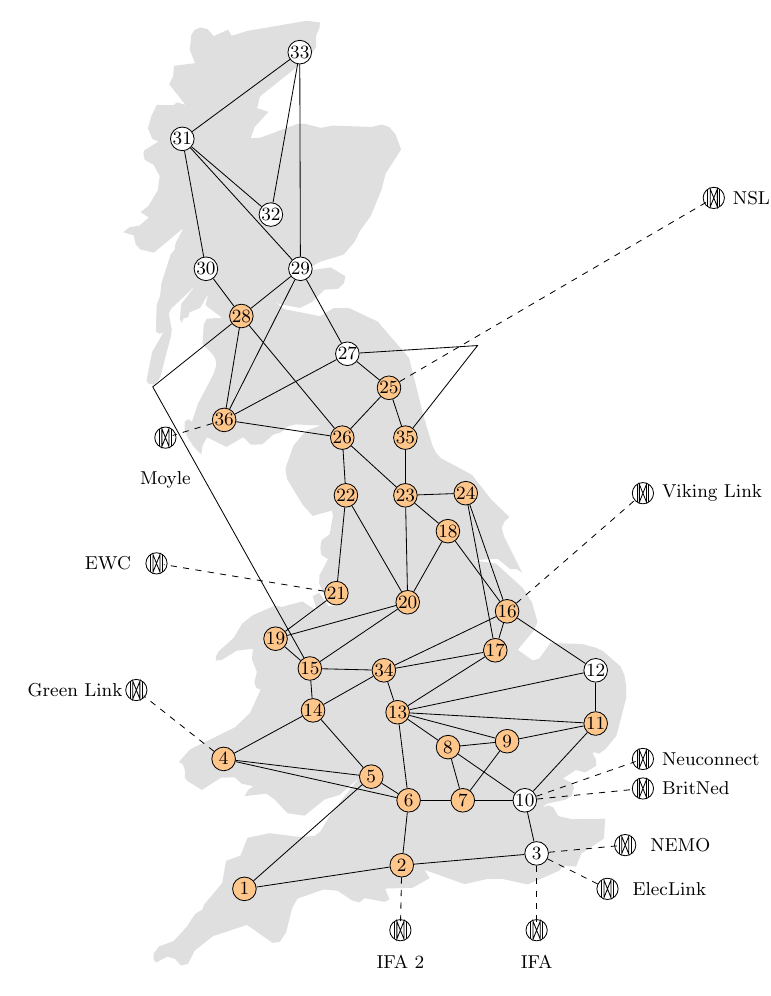}
		\label{fig:network_storage_only}
	}
	
	\caption{Comparison of the spatial expansion decisions across Cases A--D,
		over the 2030--2050 planning horizon.}
	\label{fig:network_comparison}
\end{figure*}

Figure~\ref{fig:network_comparison} compares transmission reinforcement and storage deployment decisions across Cases A--D over the 2030--2050 planning horizon. The expansion quantities are summarized in Table~\ref{tab:expansion_decisions}.

\begin{table}[H]
	\centering
	\caption{Transmission and storage expansion decisions across Cases A--D.}
	\label{tab:expansion_decisions}
	\begingroup
	\renewcommand{\arraystretch}{1.15}
	\setlength{\tabcolsep}{6pt}
	
	\begin{tabular}{@{}llrrl@{}}
		\hline
		\textbf{Case}
		& \textbf{Asset}
		& \textbf{No. built}
		& \textbf{Capacity}
		& \textbf{Unit} \\
		\hline
		Case A & Transmission & 45 & 69,150 & MW  \\
		& Storage      & 32 & 38,400 & MWh \\
		\hline
		Case B & Transmission & 45 & 69,150 & MW  \\
		& Storage      & 32 & 38,400 & MWh \\
		\hline
		Case C & Transmission & 46 & 70,850 & MW  \\
		& Storage      & 0  & 0      & MWh \\
		\hline
		Case D & Transmission & 0  & 0      & MW  \\
		& Storage      & 26 & 31,200 & MWh \\
		\hline
	\end{tabular}
	
	\endgroup
\end{table}

Case A selects 45 transmission reinforcements, providing 69,150 MW of additional capacity, and 32 storage facilities with a combined energy capacity of 38,400 MWh. Case B selects exactly the same transmission assets and storage
sites, with identical final capacities. Consequently, their spatial investment
patterns aggregated over the full planning horizon are identical. The effect
of removing the transmission preparation period is expressed through investment
timing and system costs rather than a different final asset portfolio in these
simulations. Investment decision years differ for 19 transmission assets and
eight storage facilities between Cases A and B.

Case C selects 46 transmission reinforcements with a combined additional capacity of 70,850 MW and no new storage. Its selected transmission portfolio contains all 45 assets selected in Case A plus one additional reinforcement, L40-1823. Coordinated planning therefore reduces the number of transmission reinforcements by 2.17\% and the additional transmission capacity by 1,700 MW, or 2.40\%, relative to transmission-only planning. The reduction in expansion scale is modest, while the associated improvement in system adequacy is substantial, as discussed below.

Case D selects no transmission reinforcements and 26 storage facilities
totalling 31,200 MWh. Its storage capacity is lower than the 38,400 MWh
selected in Cases A and B. Despite this deployment, Case D experiences
substantial load shedding, indicating that the selected storage-only portfolio cannot compensate for the absence of additional inter-regional transmission capacity under the modeled conditions.

\clearpage

\subsection{Comparison of Planning Strategies under Four FES Scenarios}
\label{subsec:strategy_comparison}
This subsection compares the expected costs, system adequacy, and investment decisions of Cases A--D under the four FES 2025 pathways. Table~\ref{tab:case_comparison} provides an overall quantitative comparison of the four planning strategies, including total expected system cost, interconnector cost, generation cost, renewable curtailment cost, transmission and storage investment expenditure, storage operating cost, load-shedding cost, cumulative expected load shedding, and the final MIPGap. Table~\ref{tab:expansion_decisions} complements this comparison by summarizing the corresponding transmission and storage expansion decisions, while Figures~\ref{fig:Expected_cost breakdown_across_Cases}--\ref{fig:cumulative_expansion} illustrate the cost composition, annual expected load shedding, and investment patterns in greater detail.

\begin{table}[H]\footnotesize
	\centering
	\caption{Overall comparison of Cases A--D under the four FES 2025 pathways.}
	\label{tab:case_comparison}
	\begingroup
	\renewcommand{\arraystretch}{1.15}
	\setlength{\tabcolsep}{3.5pt}
	
	\begin{tabular}{@{}lrrrrrrrrrrr@{}}
		\hline
		\textbf{Case}
		& \textbf{Total}
		& \textbf{IC}
		& \textbf{Gen.}
		& \textbf{Wind}
		& \textbf{PV}
		& \textbf{TEP}
		& \textbf{ESS}
		& \textbf{ESS}
		& \textbf{LS}
		& \textbf{LS}
		& \textbf{Gap} \\
		
		&
		& \textbf{cost}
		& \textbf{cost}
		& \textbf{curt.}
		& \textbf{curt.}
		& \textbf{inv.}
		& \textbf{inv.}
		& \textbf{OPEX}
		& \textbf{cost}
		& \textbf{(TWh)}
		& \textbf{(\%)} \\
		
		\hline
		Case A
		& 578.15 & 23.74 & 473.59 & 50.49 & 2.21
		& 4.43 & 5.79 & 7.87 & 10.03 & 1.003 & 0.25 \\
		
		Case B
		& 572.63 & 23.66 & 470.23 & 48.28 & 2.04
		& 4.37 & 5.93 & 8.09 & 10.03 & 1.003 & 0.23 \\
		
		Case C
		& 903.15 & 21.09 & 481.85 & 56.62 & 2.43
		& 4.63 & 0 & 0 & 336.54 & 33.654 & 0.29 \\
		
		Case D
		& 3163.16 & 21.22 & 551.59 & 114.31 & 6.81
		& 0 & 6.77 & 14.11 & 2448.36 & 244.836 & 0.14 \\
		
		\hline
	\end{tabular}
	
	\par
	\vspace{3pt}
	
	\begin{minipage}{0.98\textwidth}
		\footnotesize
		All cost components are in \pounds\,billion.
		IC: interconnector; Gen.: generation; curt.: curtailment;
		TEP: transmission expansion planning; ESS: energy storage system;
		LS: load shedding; Gap: final MIPGap.
	\end{minipage}
	
	\endgroup
\end{table}

As shown in Table~\ref{tab:case_comparison}, Cases A and B produce very similar overall system outcomes. Their total expected costs are \pounds578.15 billion and \pounds572.63 billion, respectively, while both cases incur the same cumulative expected load shedding of 1.003 TWh and the same load-shedding penalty of \pounds10.03 billion. The similarity in these results is consistent with their identical final transmission and storage portfolios, indicating that the main effect of removing the transmission preparation period is on the timing of investment and the associated discounted cost rather than on the final
expansion scale.

Case C, which permits transmission expansion but no additional storage, has a substantially higher total expected cost of \pounds903.15 billion. Although its generation, wind-curtailment, PV-curtailment, and transmission investment costs remain of a similar order to those of Cases A and B, its cumulative expected load shedding increases to 33.654 TWh, corresponding to a load-shedding cost of \pounds336.54 billion. This contrast highlights that the principal difference between coordinated and transmission-only planning is not the magnitude of transmission investment expenditure itself, but the ability of storage to improve system adequacy and reduce unserved energy.

Case D shows the most severe deterioration in system performance.
Without transmission expansion, cumulative expected load shedding rises to 244.836 TWh and the associated load-shedding cost reaches
\pounds2,448.36 billion. Generation cost also increases to \pounds551.59 billion, while wind- and PV-curtailment costs rise to \pounds114.31 billion and \pounds6.81 billion, respectively. Although Case D includes \pounds6.77 billion of storage investment and \pounds14.11 billion of storage operating expenditure, these additional storage resources are not sufficient to compensate for the absence of transmission reinforcement under the modeled conditions.

The final MIPGap values range from 0.14\% to 0.29\% across all four cases, indicating that the reported solutions are consistently close to optimal within the specified solver tolerance. Overall, Table~\ref{tab:case_comparison} shows that the dominant differences among the planning strategies arise from system adequacy and load-shedding performance, whereas most non-load-shedding cost components remain comparatively similar across Cases A--C.

\subsubsection{Expected Cost Comparison}
\begin{figure}[H]
	\centering

    \begin{tikzpicture}
        \begin{axis}[
            width=\textwidth,
            hide axis,
            xmin=0, xmax=1,
            ymin=0, ymax=1,
            legend columns=3,
            legend style={
                draw=none,
                font=\small,
                /tikz/every even column/.append style={
                    column sep=0.8cm
                }
            }
        ]

            \addlegendimage{area legend, fill=blue!70!black, draw=none}
            \addlegendentry{Generator}

            \addlegendimage{area legend, fill=orange!80, draw=none}
            \addlegendentry{Interconnector}

            \addlegendimage{area legend, fill=green!60!black, draw=none}
            \addlegendentry{Wind curtailment}

            \addlegendimage{area legend, fill=red!60, draw=none}
            \addlegendentry{PV curtailment}

            \addlegendimage{area legend, fill=purple!70, draw=none}
            \addlegendentry{Line expansion}

            \addlegendimage{area legend, fill=brown!70, draw=none}
            \addlegendentry{Storage operation}

            \addlegendimage{area legend, fill=magenta!60, draw=none}
            \addlegendentry{Storage investment}

            \addlegendimage{area legend, fill=gray!55, draw=none}
            \addlegendentry{Load shedding}

        \end{axis}
    \end{tikzpicture}
	\subfloat[Costs with load shedding.]{
		\resizebox{0.47\linewidth}{!}{
			\begin{tikzpicture}
            \begin{axis}[
                ybar stacked,
                width=0.49\textwidth,
                height=7cm,
                bar width=0.55cm,
                symbolic x coords={Case A,Case B,Case C,Case D},
                xtick=data,
                ymin=0,
                ymax=3500,
                ylabel={Expected cost (£ billion)},
                ylabel style={font=\small},
                xlabel={},
                tick label style={font=\small},
                legend style={
                    draw=none,
                    font=\small
                },
                enlarge x limits=0.15,
                axis line style={black},
                ymajorgrids=true,
                grid style={dashed},
            ]

                \addplot+[fill=blue!70!black, draw=none]
                    coordinates {
                        (Case A,473.588095)
                        (Case B,470.232397)
                        (Case C,481.849334)
                        (Case D,551.588313)
                    };

                \addplot+[fill=orange!80, draw=none]
                    coordinates {
                        (Case A,23.739070)
                        (Case B,23.663731)
                        (Case C,21.086740)
                        (Case D,21.218939)
                    };

                \addplot+[fill=green!60!black, draw=none]
                    coordinates {
                        (Case A,50.491610)
                        (Case B,48.280346)
                        (Case C,56.624666)
                        (Case D,114.306034)
                    };

                \addplot+[fill=red!60, draw=none]
                    coordinates {
                        (Case A,2.209769)
                        (Case B,2.038133)
                        (Case C,2.426597)
                        (Case D,6.810757)
                    };

                \addplot+[fill=purple!70, draw=none]
                    coordinates {
                        (Case A,4.431547)
                        (Case B,4.370239)
                        (Case C,4.630150)
                        (Case D,0)
                    };

                \addplot+[fill=brown!70, draw=none]
                    coordinates {
                        (Case A,7.871309)
                        (Case B,8.085161)
                        (Case C,0)
                        (Case D,14.111383)
                    };

                \addplot+[fill=magenta!60, draw=none]
                    coordinates {
                        (Case A,5.792192)
                        (Case B,5.929285)
                        (Case C,0)
                        (Case D,6.769106)
                    };

                \addplot+[fill=gray!55, draw=none]
                    coordinates {
                        (Case A,10.025941)
                        (Case B,10.025941)
                        (Case C,336.535023)
                        (Case D,2448.360225)
                    };

            \end{axis}
        \end{tikzpicture}
		}
		\label{fig:costs_including_the load-shedding_penalty}
	}
	\hfill
	\subfloat[Costs without load shedding.]{
		\resizebox{0.47\linewidth}{!}{
			\begin{tikzpicture}
            \begin{axis}[
                ybar stacked,
                width=0.49\textwidth,
                height=7cm,
                bar width=0.55cm,
                symbolic x coords={Case A,Case B,Case C,Case D},
                xtick=data,
                ymin=0,
                ymax=800,
                ylabel={Expected cost (£ billion)},
                ylabel style={font=\small},
                xlabel={},
                tick label style={font=\small},
                enlarge x limits=0.15,
                axis line style={black},
                ymajorgrids=true,
                grid style={dashed},
            ]

                \addplot+[fill=blue!70!black, draw=none]
                    coordinates {
                        (Case A,473.588095)
                        (Case B,470.232397)
                        (Case C,481.849334)
                        (Case D,551.588313)
                    };

                \addplot+[fill=orange!80, draw=none]
                    coordinates {
                        (Case A,23.739070)
                        (Case B,23.663731)
                        (Case C,21.086740)
                        (Case D,21.218939)
                    };

                \addplot+[fill=green!60!black, draw=none]
                    coordinates {
                        (Case A,50.491610)
                        (Case B,48.280346)
                        (Case C,56.624666)
                        (Case D,114.306034)
                    };

                \addplot+[fill=red!60, draw=none]
                    coordinates {
                        (Case A,2.209769)
                        (Case B,2.038133)
                        (Case C,2.426597)
                        (Case D,6.810757)
                    };

                \addplot+[fill=purple!70, draw=none]
                    coordinates {
                        (Case A,4.431547)
                        (Case B,4.370239)
                        (Case C,4.630150)
                        (Case D,0)
                    };

                \addplot+[fill=brown!70, draw=none]
                    coordinates {
                        (Case A,7.871309)
                        (Case B,8.085161)
                        (Case C,0)
                        (Case D,14.111383)
                    };

                \addplot+[fill=magenta!60, draw=none]
                    coordinates {
                        (Case A,5.792192)
                        (Case B,5.929285)
                        (Case C,0)
                        (Case D,6.769106)
                    };

            \end{axis}
        \end{tikzpicture}
		}
		\label{fig:costs_excluding_the_load-shedding_penalty}
	}
	
	\caption{Expected cost breakdown across Cases A--D:
		(a) costs including the load-shedding penalty;
		(b) costs excluding the load-shedding penalty.}
	\label{fig:Expected_cost breakdown_across_Cases}
\end{figure}
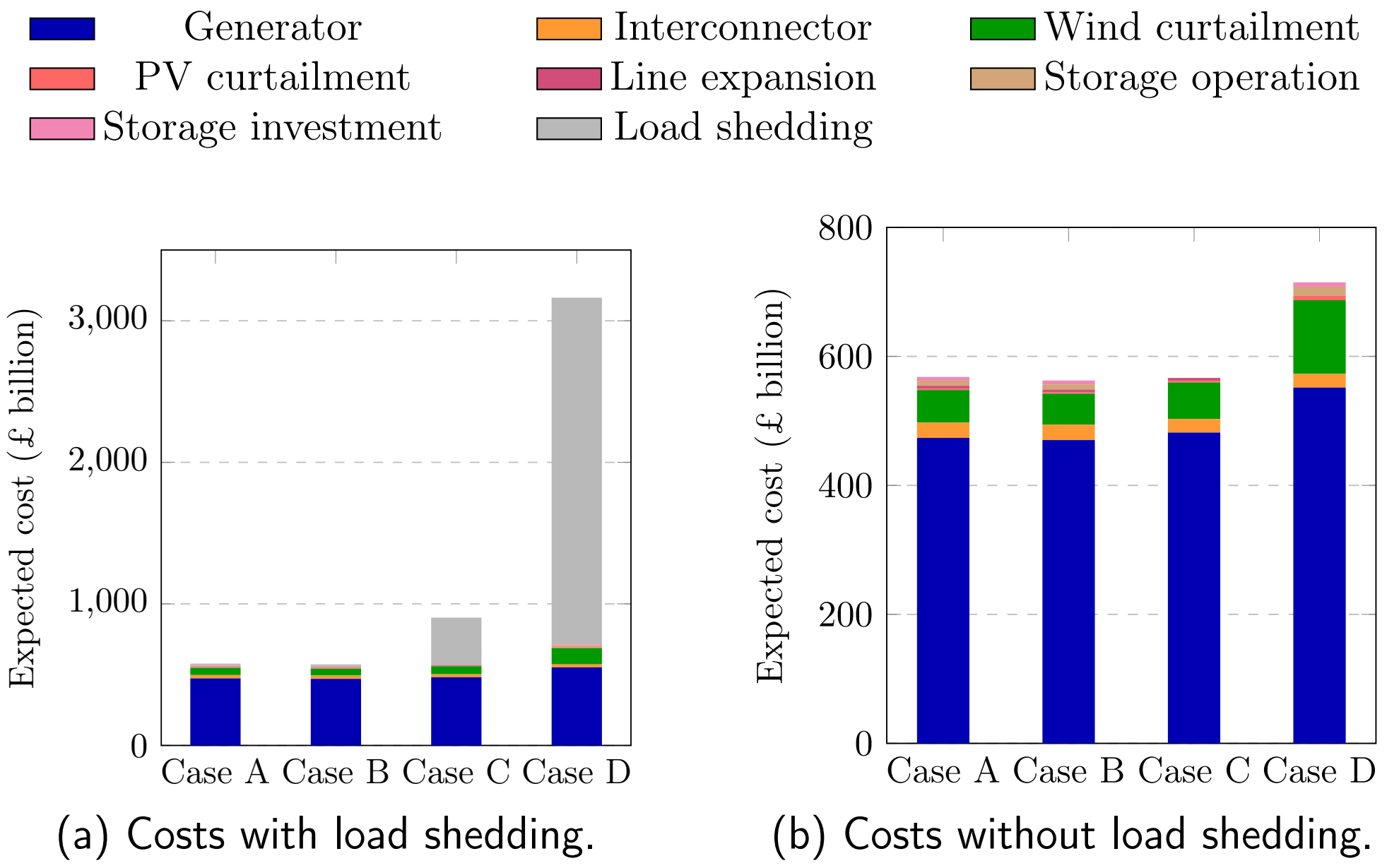

Figure~\ref{fig:Expected_cost breakdown_across_Cases}\subref{fig:costs_including_the load-shedding_penalty}
shows that the large difference between coordinated planning and
transmission-only planning is dominated by the load-shedding component.
Relative to Case C, Case A reduces total expected cost by
\pounds325.00 billion, or 35.99\%. The load-shedding penalty falls from
\pounds336.54 billion in Case C to \pounds10.03 billion in Case A.
Generation, wind-curtailment, and PV-curtailment costs are also lower in
Case A, although these differences are much smaller than the reduction in
the load-shedding penalty.

To isolate the remaining cost components,
Figure~\ref{fig:Expected_cost breakdown_across_Cases}
\subref{fig:costs_excluding_the_load-shedding_penalty}
shows the same cost breakdown after excluding load shedding. In this panel,
the non-load-shedding costs of Cases A--C are much closer than their total
costs in Figure~\ref{fig:Expected_cost breakdown_across_Cases}
\subref{fig:costs_including_the load-shedding_penalty}
would suggest. After excluding load shedding, the total costs are
\pounds568.12 billion for Case A and \pounds566.62 billion for Case C.
This indicates that the overall economic benefit of coordinated planning is
driven primarily by improved supply adequacy rather than by a substantial
reduction in the combined generation, curtailment, transmission-investment,
and storage-related cost components.

Case D exhibits a markedly different cost structure. As illustrated in
Figure~\ref{fig:Expected_cost breakdown_across_Cases}
\subref{fig:costs_including_the load-shedding_penalty},
its load-shedding penalty reaches \pounds2,448.36 billion, accounting for
approximately 77.40\% of its total expected cost. This component dominates
the stacked cost profile and explains why the total cost of Case D is much
higher than those of the other planning strategies. When load shedding is
excluded in
Figure~\ref{fig:Expected_cost breakdown_across_Cases}
\subref{fig:costs_excluding_the_load-shedding_penalty},
the remaining cost components become easier to compare, and Case D still
shows higher generation and renewable-curtailment costs than the coordinated
cases. This suggests that storage-only expansion cannot fully substitute for
additional inter-regional transmission capacity under the modelled conditions.

Interconnector net costs are positive in all four cases and therefore
contribute positively to the stacked cost totals in both panels of
Figure~\ref{fig:Expected_cost breakdown_across_Cases}. Overall,
Figure~\ref{fig:Expected_cost breakdown_across_Cases}
demonstrates that the principal economic advantage of coordinated
transmission and storage planning arises from its ability to reduce unserved
energy, while differences in the remaining cost components are comparatively
modest.

\subsubsection{System Adequacy}
\begin{figure}[H]
\subfloat[Load shedding volume in Cases A - D.]{\begin{tikzpicture}
\begin{axis}[
    width=0.95\textwidth,
    height=3.4cm,
    xlabel={Year},
    ylabel={Annual LS (TWh)},
    xmin=2030, xmax=2050,
    ymin=0, ymax=25,
    xtick={2030,2032,2034,2036,2038,2040,2042,2044,2046,2048,2050},
    ytick={0,5,10,15,20,25},
    grid=both,
    grid style={dashed, gray!30},
    axis lines=left,
    legend columns=4,
    legend style={
        at={(0.5,1.05)},
        anchor=south,
        draw=none,
        font=\small
    },
    legend cell align={left},
    tick label style={font=\small},
    label style={font=\small},
    /pgf/number format/.cd, 1000 sep={},
    scaled ticks=false,
]

\addplot[
    color=blue!70!black,
    mark=o,
    mark options={fill=blue!70!black, scale=0.7},
    line width=0.8pt
] coordinates {
(2030,0.000000) (2031,0.000000) (2032,0.000000) (2033,0.000000) (2034,0.000000)
(2035,0.000000) (2036,0.000000) (2037,0.000000) (2038,0.000000) (2039,0.000000)
(2040,0.000000) (2041,0.000000) (2042,0.000000) (2043,0.000000) (2044,0.000000)
(2045,0.278145) (2046,0.460909) (2047,0.263540) (2048,0.000000) (2049,0.000000)
(2050,0.000000)
};
\addlegendentry{Case A}

\addplot[
    color=orange,
    mark=diamond,
    mark options={fill=orange, scale=0.7},
    line width=0.8pt,
    densely dashed
] coordinates {
(2030,0.000000) (2031,0.000000) (2032,0.000000) (2033,0.000000) (2034,0.000000)
(2035,0.000000) (2036,0.000000) (2037,0.000000) (2038,0.000000) (2039,0.000000)
(2040,0.000000) (2041,0.000000) (2042,0.000000) (2043,0.000000) (2044,0.000000)
(2045,0.278145) (2046,0.460909) (2047,0.263540) (2048,0.000000) (2049,0.000000)
(2050,0.000000)
};
\addlegendentry{Case B}

\addplot[
    color=teal,
    mark=triangle*,
    mark options={fill=teal, scale=0.9},
    line width=0.9pt
] coordinates {
(2030,0.000000) (2031,0.000000) (2032,0.000000) (2033,0.000000) (2034,0.000000)
(2035,0.000000) (2036,0.000000) (2037,0.082551) (2038,0.795355) (2039,2.320432)
(2040,3.980466) (2041,0.988541) (2042,1.765883) (2043,2.541271) (2044,2.995891)
(2045,3.443806) (2046,3.666220) (2047,3.577801) (2048,3.161241) (2049,2.615048)
(2050,1.718997)
};
\addlegendentry{Case C}

\addplot[
    color=red,
    mark=*,
    mark options={fill=red, scale=0.7},
    line width=0.9pt
] coordinates {
(2030,0.000000) (2031,0.000000) (2032,0.000000) (2033,0.000000) (2034,0.000000)
(2035,1.371143) (2036,4.060307) (2037,7.099815) (2038,10.386484) (2039,14.347065)
(2040,19.410585) (2041,11.476658) (2042,14.345280) (2043,17.593919) (2044,19.698191)
(2045,21.103993) (2046,21.946782) (2047,22.286924) (2048,21.455316) (2049,20.356660)
(2050,17.896902)
};
\addlegendentry{Case D}

\end{axis}
\end{tikzpicture}}\label{fig:load_shedding_abcd}\\
\subfloat[Load shedding volume in Cases A - C.]{\begin{tikzpicture}
\begin{axis}[
    width=0.95\textwidth,
    height=3.4cm,
    xlabel={Year},
    ylabel={Annual LS (TWh)},
    xmin=2030, xmax=2050,
    ymin=0, ymax=6,
    xtick={2030,2032,2034,2036,2038,2040,2042,2044,2046,2048,2050},
    ytick={0,1,2,3,4,5},
    grid=both,
    grid style={dashed, gray!30},
    axis lines=left,
    legend columns=4,
    legend style={
        at={(0.5,1.05)},
        anchor=south,
        draw=none,
        font=\small
    },
    legend cell align={left},
    tick label style={font=\small},
    label style={font=\small},
    /pgf/number format/.cd, 1000 sep={},
    scaled ticks=false,
]

\addplot[
    color=blue!70!black,
    mark=o,
    mark options={fill=blue!70!black, scale=0.7},
    line width=0.8pt
] coordinates {
(2030,0.000000) (2031,0.000000) (2032,0.000000) (2033,0.000000) (2034,0.000000)
(2035,0.000000) (2036,0.000000) (2037,0.000000) (2038,0.000000) (2039,0.000000)
(2040,0.000000) (2041,0.000000) (2042,0.000000) (2043,0.000000) (2044,0.000000)
(2045,0.278145) (2046,0.460909) (2047,0.263540) (2048,0.000000) (2049,0.000000)
(2050,0.000000)
};

\addplot[
    color=orange,
    mark=diamond,
    mark options={fill=orange, scale=0.7},
    line width=0.8pt,
    densely dashed
] coordinates {
(2030,0.000000) (2031,0.000000) (2032,0.000000) (2033,0.000000) (2034,0.000000)
(2035,0.000000) (2036,0.000000) (2037,0.000000) (2038,0.000000) (2039,0.000000)
(2040,0.000000) (2041,0.000000) (2042,0.000000) (2043,0.000000) (2044,0.000000)
(2045,0.278145) (2046,0.460909) (2047,0.263540) (2048,0.000000) (2049,0.000000)
(2050,0.000000)
};

\addplot[
    color=teal,
    mark=triangle*,
    mark options={fill=teal, scale=0.9},
    line width=0.9pt
] coordinates {
(2030,0.000000) (2031,0.000000) (2032,0.000000) (2033,0.000000) (2034,0.000000)
(2035,0.000000) (2036,0.000000) (2037,0.082551) (2038,0.795355) (2039,2.320432)
(2040,3.980466) (2041,0.988541) (2042,1.765883) (2043,2.541271) (2044,2.995891)
(2045,3.443806) (2046,3.666220) (2047,3.577801) (2048,3.161241) (2049,2.615048)
(2050,1.718997)
};

\end{axis}
\end{tikzpicture}}\label{fig:load_shedding_abc_zoomed}
	\caption{Annual expected load shedding over the 2030--2050 planning horizon: (a) Cases A--D; (b) Cases A--C shown with a reduced vertical range. The profiles of Cases A and B overlap within numerical precision.}
\label{fig:annual_load_shedding_cases}
\end{figure}
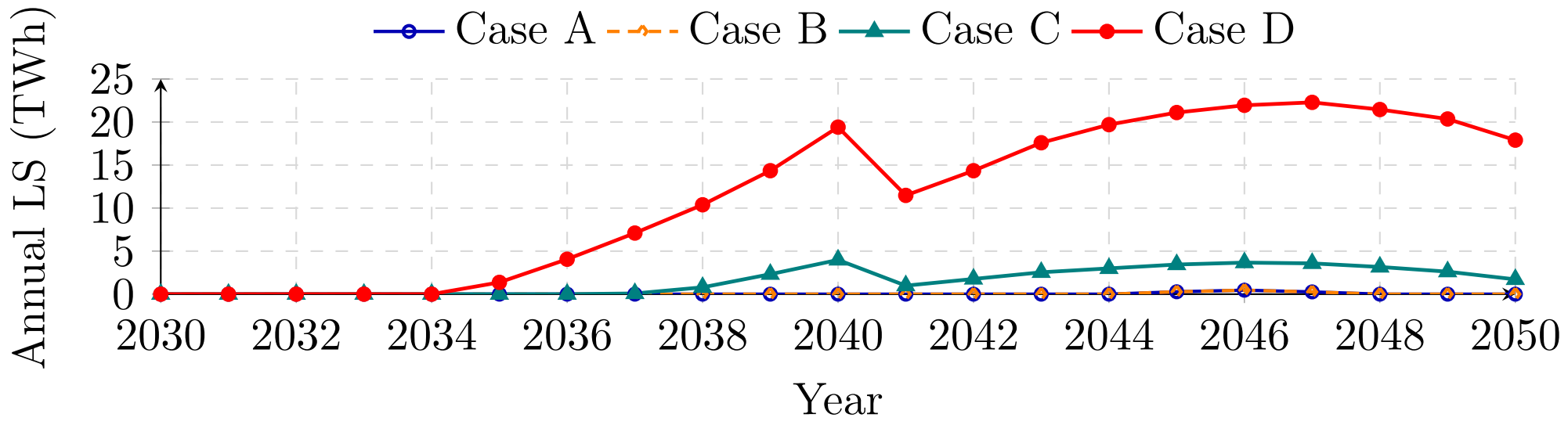
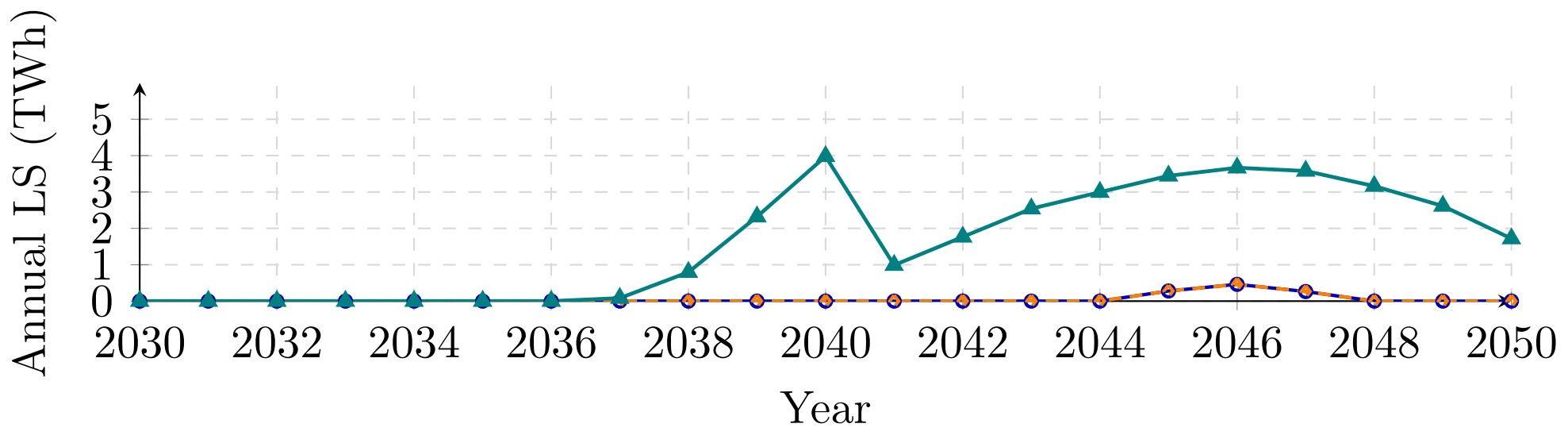

Cases A and B have effectively identical annual expected load shedding profiles,
with a cumulative value of 1.003 TWh in each case. Material load shedding occurs
in 2045--2047, reaching approximately 0.278, 0.461, and 0.264 TWh, respectively.
Their overlapping curves in Figure~\ref{fig:annual_load_shedding_cases} are distinguished by different line styles
and markers. Coordinated planning substantially reduces load shedding but does
not eliminate it under the modelled conditions.

Case C incurs cumulative expected load shedding of 33.654 TWh despite adding
70.850 GW of transmission capacity. Material load shedding begins in 2037 and
peaks at approximately 3.980 TWh in 2040. Compared with Case C, Case A reduces
cumulative expected load shedding by approximately 97.02\%, illustrating the
adequacy benefit associated with including storage in the coordinated portfolio.

Case D incurs cumulative expected load shedding of 244.836 TWh despite
deploying 31.2 GWh of storage. Material load shedding begins in 2035 and
reaches its maximum of approximately 22.287 TWh in 2047. The result indicates
that storage-only expansion cannot maintain the same level of adequacy as
coordinated transmission and storage planning in this system.

\subsubsection{Investment Scale and Timing}
Cases A, B, and C select 45, 45, and 46 transmission reinforcements, providing
69.150, 69.150, and 70.850 GW of additional capacity, respectively. Cases A
and B each select 32 storage facilities totalling 38.4 GWh, whereas Case D
selects 26 facilities totalling 31.2 GWh. Relative to Case C, Case A reduces
transmission reinforcement count and capacity by 2.17\% and 2.40\%, respectively.

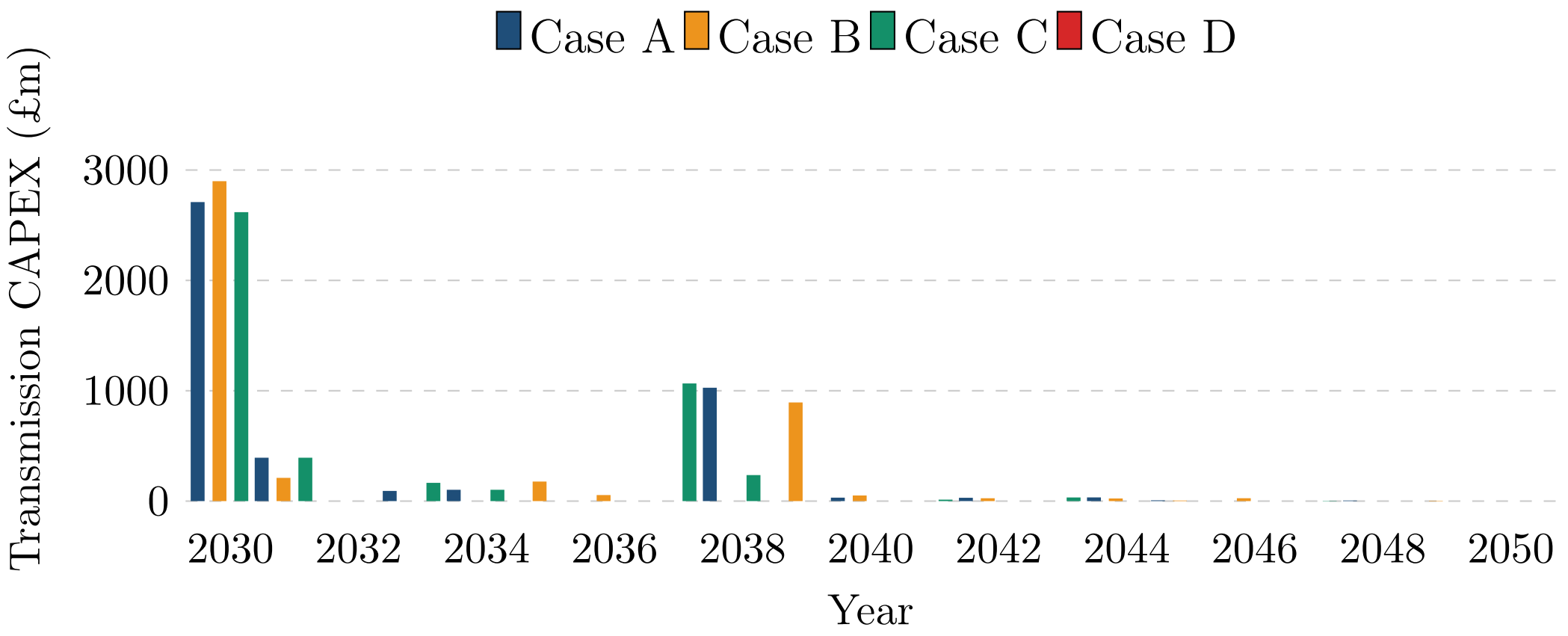
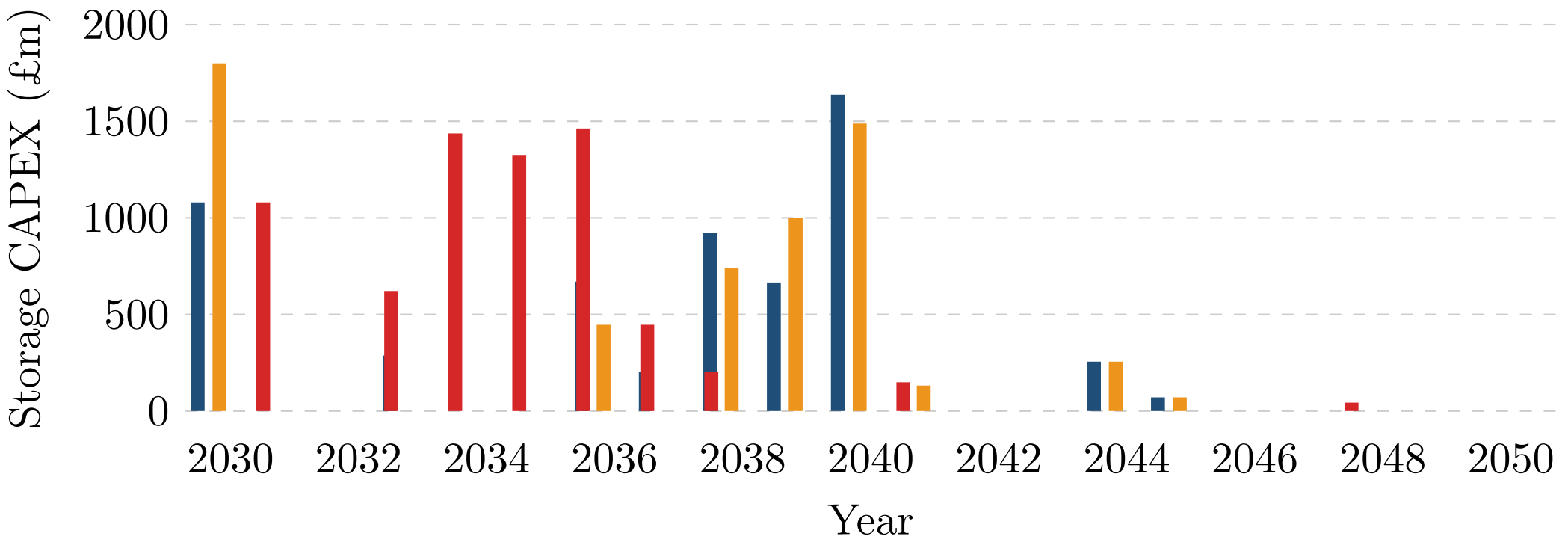
\begin{figure}[H]
\subfloat[transmission investment cost.]{\definecolor{caseA}{RGB}{31,78,121}
\definecolor{caseB}{RGB}{237,148,29}
\definecolor{caseC}{RGB}{20,144,105}
\definecolor{caseD}{RGB}{214,39,40}

\begin{tikzpicture}
\begin{axis}[
    ybar, bar width=3.5pt,
    width=\textwidth, height=5cm,
    xtick={2030,2032,2034,2036,2038,2040,2042,2044,2046,2048,2050},
    xlabel={Year},
    ylabel={Transmission CAPEX (\pounds m)},
    ymin=0, ymax=3500, ytick={0,1000,...,3000},
    enlarge x limits=0.035,
    ymajorgrids, grid style={dashed, gray!40},
    axis line style={draw=none}, tick style={draw=none},
    legend style={at={(0.5,1.1)}, anchor=south, legend columns=4, draw=none},
    tick label style={font=\small},
    label style={font=\small},
    /pgf/number format/.cd, 1000 sep={},
    scaled ticks=false,
]
\addplot[fill=caseA, draw=none,single ybar legend] coordinates {
(2030,2709.679136) (2031,393.065066) (2032,0) (2033,92.207223) (2034,102.536704)
(2035,0) (2036,0) (2037,0) (2038,1027.212629) (2039,0)
(2040,31.246421) (2041,0) (2042,30.018034) (2043,0) (2044,33.415186)
(2045,6.883988) (2046,0) (2047,0) (2048,5.282315) (2049,0) (2050,0)};
\addplot[fill=caseB, draw=none,single ybar legend] coordinates {
(2030,2898.980221) (2031,210.520379) (2032,0) (2033,0) (2034,0)
(2035,177.084131) (2036,54.678572) (2037,0) (2038,0) (2039,893.778249)
(2040,50.928923) (2041,0) (2042,25.815680) (2043,0) (2044,23.790496)
(2045,5.133883) (2046,26.066050) (2047,0) (2048,0) (2049,3.462121) (2050,0)};
\addplot[fill=caseC, draw=none,single ybar legend] coordinates {
(2030,2618.287695) (2031,393.065066) (2032,0) (2033,165.409247) (2034,102.536704)
(2035,0) (2036,0) (2037,1066.164884) (2038,235.751811) (2039,0)
(2040,0) (2041,13.219039) (2042,0) (2043,32.904052) (2044,0)
(2045,0) (2046,0) (2047,2.811453) (2048,0) (2049,0) (2050,0)};
\addplot[fill=caseD, draw=none,single ybar legend] coordinates {
(2030,0) (2031,0) (2032,0) (2033,0) (2034,0)
(2035,0) (2036,0) (2037,0) (2038,0) (2039,0)
(2040,0) (2041,0) (2042,0) (2043,0) (2044,0)
(2045,0) (2046,0) (2047,0) (2048,0) (2049,0) (2050,0)};
\legend{Case A, Case B, Case C, Case D}
\end{axis}
\end{tikzpicture}}\\
\subfloat[storage investment cost.]{\definecolor{caseA}{RGB}{31,78,121}
\definecolor{caseB}{RGB}{237,148,29}
\definecolor{caseC}{RGB}{20,144,105}
\definecolor{caseD}{RGB}{214,39,40}

\begin{tikzpicture}
\begin{axis}[
    ybar, bar width=3.5pt,
    width=\textwidth, height=5cm,
    xtick={2030,2032,2034,2036,2038,2040,2042,2044,2046,2048,2050},
    xlabel={Year},
    ylabel={Storage CAPEX (\pounds m)},
    ymin=0, ymax=2000, ytick={0,500,...,2000},
    enlarge x limits=0.035,
    ymajorgrids, grid style={dashed, gray!40},
    axis line style={draw=none}, tick style={draw=none},
    legend style={at={(0.5,1.1)}, anchor=south, legend columns=4, draw=none},
    tick label style={font=\small},
    label style={font=\small},
    /pgf/number format/.cd, 1000 sep={},
    scaled ticks=false,
]
\addplot[fill=caseA, draw=none,single ybar legend] coordinates {
(2030,1080.000000) (2031,0) (2032,0) (2033,287.427996) (2034,0)
(2035,0) (2036,669.725542) (2037,203.515020) (2038,922.797945) (2039,665.362366)
(2040,1637.072579) (2041,0) (2042,0) (2043,0) (2044,255.551343)
(2045,70.739437) (2046,0) (2047,0) (2048,0) (2049,0) (2050,0)};
\addplot[fill=caseB, draw=none,single ybar legend] coordinates {
(2030,1800.000000) (2031,0) (2032,0) (2033,0) (2034,0)
(2035,0) (2036,446.483695) (2037,0) (2038,738.238356) (2039,998.043548)
(2040,1488.247799) (2041,131.980536) (2042,0) (2043,0) (2044,255.551343)
(2045,70.739437) (2046,0) (2047,0) (2048,0) (2049,0) (2050,0)};
\addplot[fill=caseC, draw=none,single ybar legend] coordinates {
(2030,0) (2031,0) (2032,0) (2033,0) (2034,0)
(2035,0) (2036,0) (2037,0) (2038,0) (2039,0)
(2040,0) (2041,0) (2042,0) (2043,0) (2044,0)
(2045,0) (2046,0) (2047,0) (2048,0) (2049,0) (2050,0)};
\addplot[fill=caseD, draw=none,single ybar legend] coordinates {
(2030,1080.000000) (2031,0) (2032,621.255216) (2033,1437.139978) (2034,1325.793326)
(2035,1462.661036) (2036,446.483695) (2037,203.515020) (2038,0) (2039,0)
(2040,148.824780) (2041,0) (2042,0) (2043,0) (2044,0)
(2045,0) (2046,0) (2047,43.432921) (2048,0) (2049,0) (2050,0)};
\end{axis}
\end{tikzpicture}}
\caption{Annual effective CAPEX present values across Cases A--D:
	(a) transmission investment cost;
	(b) storage investment cost.}
\label{fig:annual_investment_expenditure}
\end{figure}

Figure~\ref{fig:annual_investment_expenditure} presents annual investment decisions weighted by their effective
CAPEX present values. Transmission investment is strongly concentrated in
2030, when Cases A, B, and C select 28, 29, and 27 reinforcements,
respectively. Further investment occurs later in the horizon, including
substantial transmission CAPEX PV in 2038, 2039, and 2037 for Cases A, B,
and C, respectively.

Storage investment follows a different temporal pattern. Cases A and B select
only three and five storage facilities in 2030, respectively, and each adds
20 facilities during 2038--2040. By contrast, Case D selects 21 of its
26 facilities by the end of 2035. Storage deployment therefore cannot be
characterized uniformly as predominantly occurring at the start of the
planning horizon.

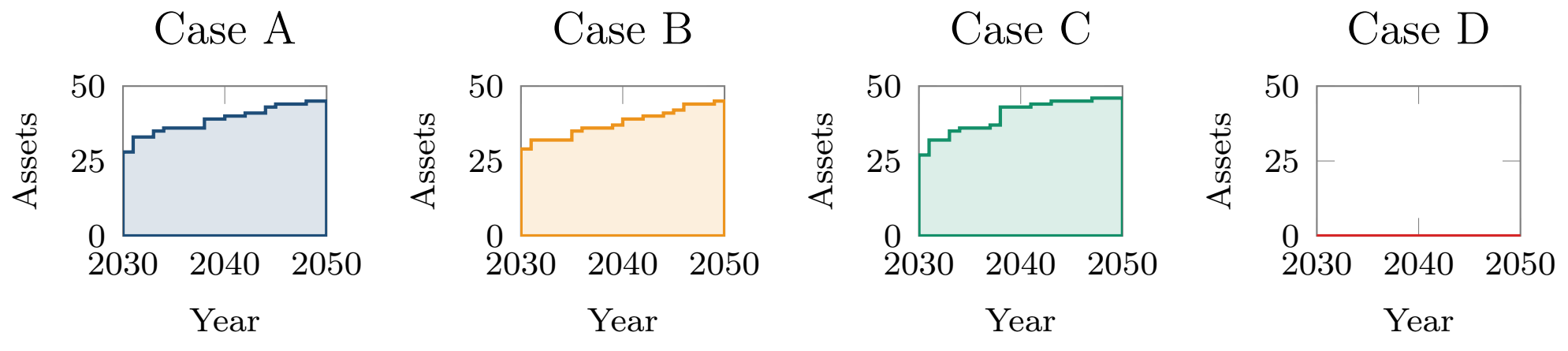
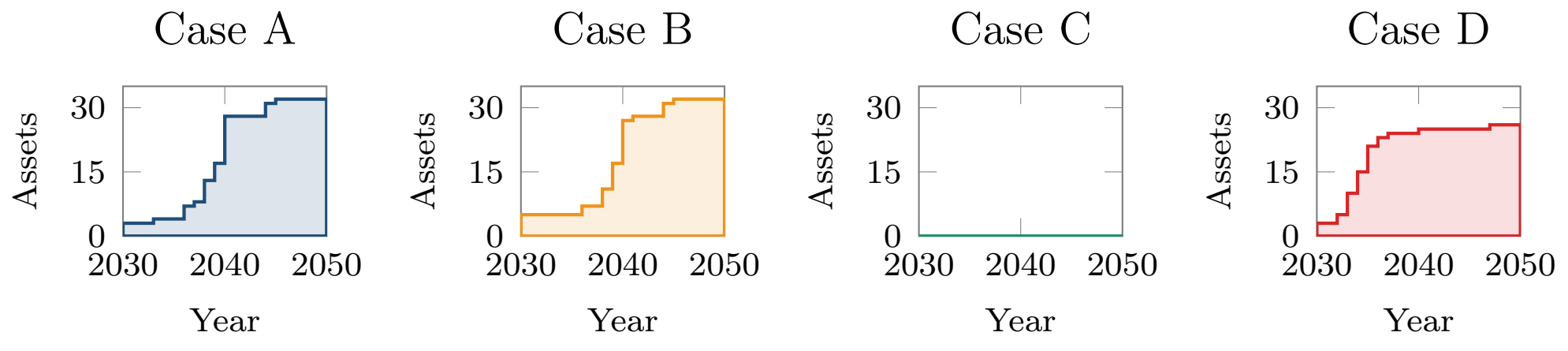
\begin{figure}[H]
\centering
\subfloat[Transmission investment decisions.]{\definecolor{caseA}{RGB}{31,78,121}
\definecolor{caseB}{RGB}{237,148,29}
\definecolor{caseC}{RGB}{20,144,105}
\definecolor{caseD}{RGB}{214,39,40}

\pgfplotsset{
    smallmult/.style={
        width=0.24\textwidth, 
        xmin=2030, xmax=2050, ymin=0, ymax=50,
        xtick={2030,2040,2050}, ytick={0,25,50},
        tick label style={font=\scriptsize},
        title style={font=\small},
        xlabel={Year}, ylabel={Assets},
        label style={font=\scriptsize},
        axis line style={gray}, tick style={gray},
    }
}

\begin{tikzpicture}
\begin{axis}[smallmult, title={Case A},
    /pgf/number format/.cd, 1000 sep={},
    scaled ticks=false,]
\addplot[thick, color=caseA, fill=caseA!15] coordinates {
(2030,0)(2030,28)(2031,28)(2031,33)(2032,33)(2032,33)(2033,33)(2033,35)
(2034,35)(2034,36)(2035,36)(2035,36)(2036,36)(2036,36)(2037,36)(2037,36)
(2038,36)(2038,39)(2039,39)(2039,39)(2040,39)(2040,40)(2041,40)(2041,40)
(2042,40)(2042,41)(2043,41)(2043,41)(2044,41)(2044,43)(2045,43)(2045,44)
(2046,44)(2046,44)(2047,44)(2047,44)(2048,44)(2048,45)(2049,45)(2049,45)
(2050,45)(2050,0)
} \closedcycle;
\end{axis}
\end{tikzpicture}
\begin{tikzpicture}
\begin{axis}[smallmult, title={Case B},
    /pgf/number format/.cd, 1000 sep={},
    scaled ticks=false,]
\addplot[thick, color=caseB, fill=caseB!15] coordinates {
(2030,0)(2030,29)(2031,29)(2031,32)(2032,32)(2032,32)(2033,32)(2033,32)
(2034,32)(2034,32)(2035,32)(2035,35)(2036,35)(2036,36)(2037,36)(2037,36)
(2038,36)(2038,36)(2039,36)(2039,37)(2040,37)(2040,39)(2041,39)(2041,39)
(2042,39)(2042,40)(2043,40)(2043,40)(2044,40)(2044,41)(2045,41)(2045,42)
(2046,42)(2046,44)(2047,44)(2047,44)(2048,44)(2048,44)(2049,44)(2049,45)
(2050,45)(2050,0)
} \closedcycle;
\end{axis}
\end{tikzpicture}
\begin{tikzpicture}
\begin{axis}[smallmult, title={Case C},
    /pgf/number format/.cd, 1000 sep={},
    scaled ticks=false,]
\addplot[thick, color=caseC, fill=caseC!15] coordinates {
(2030,0)(2030,27)(2031,27)(2031,32)(2032,32)(2032,32)(2033,32)(2033,35)
(2034,35)(2034,36)(2035,36)(2035,36)(2036,36)(2036,36)(2037,36)(2037,37)
(2038,37)(2038,43)(2039,43)(2039,43)(2040,43)(2040,43)(2041,43)(2041,44)
(2042,44)(2042,44)(2043,44)(2043,45)(2044,45)(2044,45)(2045,45)(2045,45)
(2046,45)(2046,45)(2047,45)(2047,46)(2048,46)(2048,46)(2049,46)(2049,46)
(2050,46)(2050,0)
} \closedcycle;
\end{axis}
\end{tikzpicture}
\begin{tikzpicture}
\begin{axis}[smallmult, title={Case D},
    /pgf/number format/.cd, 1000 sep={},
    scaled ticks=false,]
\addplot[thick, color=caseD, fill=caseD!15] coordinates {
(2030,0)(2050,0)
} \closedcycle;
\end{axis}
\end{tikzpicture}}\\
\centering
\subfloat[Storage investment decisions.]{\definecolor{caseA}{RGB}{31,78,121}
\definecolor{caseB}{RGB}{237,148,29}
\definecolor{caseC}{RGB}{20,144,105}
\definecolor{caseD}{RGB}{214,39,40}

\pgfplotsset{
    smallmult/.style={
        width=0.24\textwidth, 
        xmin=2030, xmax=2050, ymin=0, ymax=35,
        xtick={2030,2040,2050}, ytick={0,15,30},
        tick label style={font=\scriptsize},
        title style={font=\small},
        xlabel={Year}, ylabel={Assets},
        label style={font=\scriptsize},
        axis line style={gray}, tick style={gray},
    }
}

\begin{tikzpicture}
\begin{axis}[smallmult, title={Case A},
    /pgf/number format/.cd, 1000 sep={},
    scaled ticks=false,]
\addplot[thick, color=caseA, fill=caseA!15] coordinates {
(2030,0)(2030,3)(2031,3)(2031,3)(2032,3)(2032,3)(2033,3)(2033,4)
(2034,4)(2034,4)(2035,4)(2035,4)(2036,4)(2036,7)(2037,7)(2037,8)
(2038,8)(2038,13)(2039,13)(2039,17)(2040,17)(2040,28)(2041,28)(2041,28)
(2042,28)(2042,28)(2043,28)(2043,28)(2044,28)(2044,31)(2045,31)(2045,32)
(2046,32)(2046,32)(2047,32)(2047,32)(2048,32)(2048,32)(2049,32)(2049,32)
(2050,32)(2050,0)
} \closedcycle;
\end{axis}
\end{tikzpicture}
\begin{tikzpicture}
\begin{axis}[smallmult, title={Case B},
    /pgf/number format/.cd, 1000 sep={},
    scaled ticks=false,]
\addplot[thick, color=caseB, fill=caseB!15] coordinates {
(2030,0)(2030,5)(2031,5)(2031,5)(2032,5)(2032,5)(2033,5)(2033,5)
(2034,5)(2034,5)(2035,5)(2035,5)(2036,5)(2036,7)(2037,7)(2037,7)
(2038,7)(2038,11)(2039,11)(2039,17)(2040,17)(2040,27)(2041,27)(2041,28)
(2042,28)(2042,28)(2043,28)(2043,28)(2044,28)(2044,31)(2045,31)(2045,32)
(2046,32)(2046,32)(2047,32)(2047,32)(2048,32)(2048,32)(2049,32)(2049,32)
(2050,32)(2050,0)
} \closedcycle;
\end{axis}
\end{tikzpicture}
\begin{tikzpicture}
\begin{axis}[smallmult, title={Case C},
    /pgf/number format/.cd, 1000 sep={},
    scaled ticks=false,]
\addplot[thick, color=caseC, fill=caseC!15] coordinates {
(2030,0)(2050,0)
} \closedcycle;
\end{axis}
\end{tikzpicture}
\begin{tikzpicture}
\begin{axis}[smallmult, title={Case D},
    /pgf/number format/.cd, 1000 sep={},
    scaled ticks=false,]
\addplot[thick, color=caseD, fill=caseD!15] coordinates {
(2030,0)(2030,3)(2031,3)(2031,3)(2032,3)(2032,5)(2033,5)(2033,10)
(2034,10)(2034,15)(2035,15)(2035,21)(2036,21)(2036,23)(2037,23)(2037,24)
(2038,24)(2038,24)(2039,24)(2039,24)(2040,24)(2040,25)(2041,25)(2041,25)
(2042,25)(2042,25)(2043,25)(2043,25)(2044,25)(2044,25)(2045,25)(2045,25)
(2046,25)(2046,25)(2047,25)(2047,26)(2048,26)(2048,26)(2049,26)(2049,26)
(2050,26)(2050,0)
} \closedcycle;
\end{axis}
\end{tikzpicture}}
\caption{Cumulative numbers of selected transmission and storage investments, shown separately for each case. The upper and lower panels represent transmission and storage, respectively.}
\label{fig:cumulative_expansion}
\end{figure}

Figure~\ref{fig:cumulative_expansion} shows cumulative investment decisions separately for each case,
using common vertical scales for the same asset type to avoid overlapping
curves. Cases A and B converge to identical final asset portfolios but follow
different investment schedules. For example, the EGL1-2735 investment decision
occurs in 2038 in Case A and 2039 in Case B. These are investment decision
years, not commissioning years.

Overall, coordinated planning provides a modest reduction in transmission
expansion requirements and a substantial reduction in load shedding and total
cost relative to transmission-only planning. The comparison with Case D
further demonstrates the limitations of storage-only expansion. The Case A--B
comparison indicates that transmission preparation periods affect investment
timing and costs, even when the final portfolio remains unchanged.


\subsection{Scenario Analysis: EE versus FB}
\label{sec:scenario_analysis}
Cases E and F apply the coordinated planning framework, including transmission preparation periods, to the EE and FB
pathways, respectively. Their comparison illustrates how the selected
investment portfolio and its timing vary between these two pathways.

\subsubsection{Expected Cost Comparison}
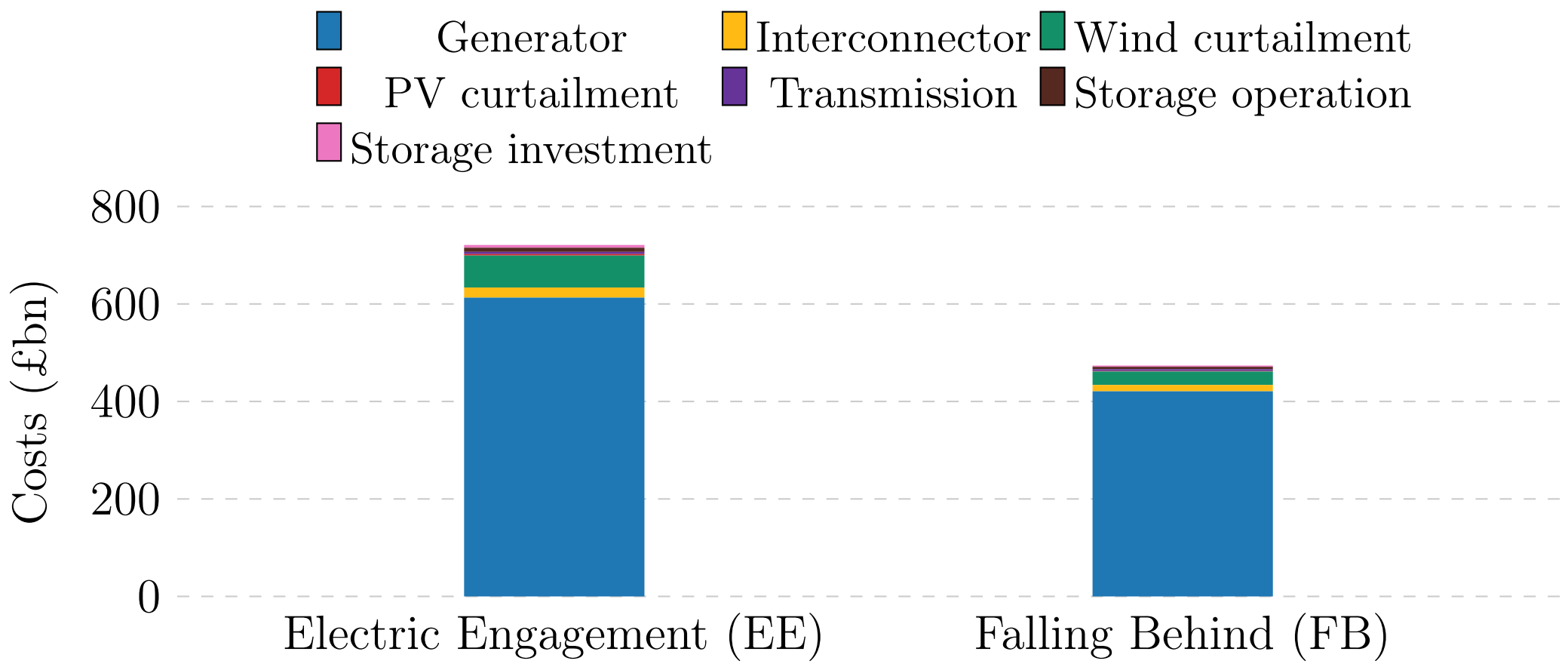
\begin{figure}[H]
	\definecolor{genCost}{RGB}{31,119,180}
\definecolor{interCost}{RGB}{255,187,20}
\definecolor{windCost}{RGB}{20,144,105}
\definecolor{pvCost}{RGB}{214,39,40}
\definecolor{transInv}{RGB}{102,51,153}
\definecolor{storOp}{RGB}{85,40,32}
\definecolor{storInv}{RGB}{237,120,193}

\begin{tikzpicture}
\begin{axis}[
    ybar stacked,
    bar width=45pt,
    width=\textwidth, height=5cm,
    symbolic x coords={Electric Engagement (EE), Falling Behind (FB)},
    xtick=data,
    ymin=0, ymax=800, ytick={0,200,...,800},
    ylabel={Costs (\pounds bn)},
    ymajorgrids, grid style={dashed, gray!40},
    axis line style={draw=none}, tick style={draw=none},
    legend style={at={(0.5,1.05)}, anchor=south, legend columns=3, draw=none, font=\small},
    enlarge x limits=0.6,
]
\addplot[fill=genCost, draw=none] coordinates {
    (Electric Engagement (EE),612.978632) (Falling Behind (FB),420.819597)};
\addplot[fill=interCost, draw=none] coordinates {
    (Electric Engagement (EE),20.799669) (Falling Behind (FB),13.031874)};
\addplot[fill=windCost, draw=none] coordinates {
    (Electric Engagement (EE),65.619019) (Falling Behind (FB),27.610243)};
\addplot[fill=pvCost, draw=none] coordinates {
    (Electric Engagement (EE),2.888125) (Falling Behind (FB),0.966042)};
\addplot[fill=transInv, draw=none] coordinates {
    (Electric Engagement (EE),4.601539) (Falling Behind (FB),3.650340)};
\addplot[fill=storOp, draw=none] coordinates {
    (Electric Engagement (EE),8.466233) (Falling Behind (FB),4.967842)};
\addplot[fill=storInv, draw=none] coordinates {
    (Electric Engagement (EE),5.664420) (Falling Behind (FB),2.547370)};
\legend{Generator, Interconnector, Wind curtailment, PV curtailment, Transmission, Storage operation, Storage investment}
\end{axis}
\end{tikzpicture}
	\caption{Expected cost breakdown under the Electric Engagement (EE) 
		and Falling Behind (FB) pathways.}
	\label{fig:ee_fb_cost}
\end{figure}

As shown in Table~\ref{tab:case_ef_comparison} and Figure~\ref{fig:ee_fb_cost}, total cost is \pounds721.018 billion
under EE and \pounds473.593 billion under FB. EE is therefore
\pounds247.424 billion, or 52.24\%, more expensive. Generation costs are
\pounds612.979 billion and \pounds420.820 billion, respectively, and their difference accounts for approximately 77.66\% of the total cost gap. Wind curtailment costs are \pounds65.619 billion under EE and
\pounds27.610 billion under FB, while PV curtailment costs are
\pounds2.888 billion and \pounds0.966 billion, respectively.

\begin{table}[t]
	\centering
	\caption{Comparison of planning results under the Electric
		Engagement and Falling Behind pathways.}
	\label{tab:case_ef_comparison}
	\begingroup
	\small
	\renewcommand{\arraystretch}{1.15}
	\setlength{\tabcolsep}{5pt}
	
	\begin{tabular}{@{}lrr@{}}
		\hline
		\textbf{Metric}
		& \textbf{Case E (EE)}
		& \textbf{Case F (FB)} \\
		\hline
		Total cost (\pounds\,billion)
		& 721.018 & 473.593 \\
		Interconnector cost (\pounds\,billion)
		& 20.800 & 13.032 \\
		Generation cost (\pounds\,billion)
		& 612.979 & 420.820 \\
		Wind curt. cost (\pounds\,billion)
		& 65.619 & 27.610 \\
		PV curt. cost (\pounds\,billion)
		& 2.888 & 0.966 \\
		TEP inv. cost (\pounds\,billion)
		& 4.602 & 3.650 \\
		Storage inv. cost (\pounds\,billion)
		& 5.664 & 2.547 \\
		Storage operation cost (\pounds\,billion)
		& 8.466 & 4.968 \\
		Load shed (MWh)
		& 0 & 0 \\
		\hline
		Transmission built
		& 43 & 40 \\
		Transmission capacity (MW)
		& 66,750 & 62,650 \\
		Storage built
		& 32 & 22 \\
		Storage capacity (MWh)
		& 38,400 & 26,400 \\
		\hline
		MIPGap (\%)
		& 0.16 & 0.33 \\
		\hline
	\end{tabular}
	
	\par\vspace{3pt}
	\begin{minipage}{\columnwidth}
		\footnotesize
		EE: Electric Engagement; FB: Falling Behind;
		curt.: curtailment; inv.: investment;
		TEP: transmission expansion planning.
		Negligible negative load shedding values due to numerical
		tolerances are reported as zero.
	\end{minipage}
	
	\endgroup
\end{table}

Combined transmission and storage investment costs are \pounds10.266 billion
under EE and \pounds6.198 billion under FB. Their absolute difference is
small relative to the total system cost gap, although storage investment
cost itself is substantially higher under EE: \pounds5.664 billion compared
with \pounds2.547 billion. Interconnector net costs are positive under both
pathways. Both cases exhibit effectively zero load shedding within numerical
tolerances under the represented operating conditions.

\subsubsection{Investment Timing}
\begin{figure}[H]
	\subfloat[]{\definecolor{eeTrans}{RGB}{31,78,151}
\definecolor{fbTrans}{RGB}{237,148,29}
\definecolor{eeStor}{RGB}{20,144,105}
\definecolor{fbStor}{RGB}{214,39,40}

\begin{tikzpicture}
\begin{axis}[
    width=\textwidth, height=5cm,
    xlabel={Year}, ylabel={CAPEX (\pounds bn)},
    xmin=2029.5, xmax=2050.5, ymin=0, ymax=3000,
    xtick={2030,2032,...,2050}, ytick={0,1000,...,3000},
yticklabels={0,1,...,3},
    x tick label style={/pgf/number format/.cd, 1000 sep=},
    tick label style={font=\small},
    label style={font=\small},
    ymajorgrids, xmajorgrids, 
    axis line style={draw=none}, tick style={draw=none},
    legend style={at={(0.5,1.05)}, anchor=south, legend columns=4, draw=none, font=\small},
    legend cell align={left},
]
\addplot[eeTrans, mark=o, thick] coordinates {
(2030,2651.423950) (2031,504.871463) (2032,0) (2033,0) (2034,102.536704)
(2035,36.233816) (2036,1168.624112) (2037,41.071354) (2038,39.021263) (2039,0)
(2040,39.128227) (2041,0) (2042,0) (2043,0) (2044,0)
(2045,0) (2046,10.112173) (2047,0) (2048,7.606705) (2049,0.909347) (2050,0)
};
\addplot[fbTrans, mark=square, dashed, thick] coordinates {
(2030,2184.485886) (2031,413.585846) (2032,107.609639) (2033,246.648522) (2034,28.769504)
(2035,87.440830) (2036,0) (2037,0) (2038,0) (2039,0)
(2040,25.375570) (2041,0) (2042,0) (2043,533.256423) (2044,0)
(2045,0) (2046,0) (2047,23.168022) (2048,0) (2049,0) (2050,0)
};
\addplot[eeStor, mark=triangle, dashdotted, thick] coordinates {
(2030,1440.000000) (2031,0) (2032,0) (2033,2299.423964) (2034,0)
(2035,0) (2036,0) (2037,0) (2038,0) (2039,0)
(2040,297.649560) (2041,527.922145) (2042,463.111158) (2043,300.563818) (2044,85.183781)
(2045,70.739437) (2046,0) (2047,173.731682) (2048,0) (2049,0) (2050,6.093985)
};
\addplot[fbStor, mark=diamond, dotted, thick] coordinates {
(2030,0) (2031,334.802499) (2032,0) (2033,0) (2034,0)
(2035,243.776839) (2036,0) (2037,203.515020) (2038,184.559589) (2039,332.681183)
(2040,0) (2041,263.961072) (2042,347.333369) (2043,100.187939) (2044,170.367562)
(2045,141.478875) (2046,113.660581) (2047,86.865841) (2048,0) (2049,18.085458) (2050,6.093985)
};
\legend{EE: transmission, FB: transmission, EE: storage, FB: storage}
\end{axis}
\end{tikzpicture}}\label{fig:ee_fb_investment_a}\\
	\subfloat[]{\definecolor{eeTrans}{RGB}{31,78,151}
\definecolor{fbTrans}{RGB}{237,148,29}
\definecolor{eeStor}{RGB}{20,144,105}
\definecolor{fbStor}{RGB}{214,39,40}

\begin{tikzpicture}
\begin{axis}[
    width=\textwidth, height=5cm,
    xlabel={Year}, ylabel={Assets},
    xmin=2029.5, xmax=2050.5, ymin=0, ymax=48,
    xtick={2030,2032,...,2050}, ytick={0,8,...,48},
    x tick label style={/pgf/number format/.cd, 1000 sep=},
    tick label style={font=\small}, label style={font=\small},
    ymajorgrids, xmajorgrids, grid style={gray!25},
    axis line style={draw=gray}, tick style={draw=gray},
    legend style={at={(0.03,0.97)}, anchor=north west, draw=gray!50, font=\small},
    legend columns=2,
]
\addplot[eeTrans, mark=o, thick, const plot] coordinates {
(2030,27) (2031,34) (2032,34) (2033,34) (2034,35) (2035,36) (2036,37)
(2037,38) (2038,39) (2039,39) (2040,40) (2041,40) (2042,40) (2043,40)
(2044,40) (2045,40) (2046,41) (2047,41) (2048,42) (2049,43) (2050,43)
};
\addplot[fbTrans, mark=square, dashed, thick, const plot] coordinates {
(2030,21) (2031,27) (2032,29) (2033,32) (2034,33) (2035,35) (2036,35)
(2037,35) (2038,35) (2039,35) (2040,36) (2041,36) (2042,36) (2043,38)
(2044,38) (2045,38) (2046,38) (2047,40) (2048,40) (2049,40) (2050,40)
};
\addplot[eeStor, mark=triangle, dashdotted, thick, const plot] coordinates {
(2030,4) (2031,4) (2032,4) (2033,12) (2034,12) (2035,12) (2036,12)
(2037,12) (2038,12) (2039,12) (2040,14) (2041,18) (2042,22) (2043,25)
(2044,26) (2045,27) (2046,27) (2047,31) (2048,31) (2049,31) (2050,32)
};
\addplot[fbStor, mark=diamond, dotted, thick, const plot] coordinates {
(2030,0) (2031,1) (2032,1) (2033,1) (2034,1) (2035,2) (2036,2)
(2037,3) (2038,4) (2039,6) (2040,6) (2041,8) (2042,11) (2043,12)
(2044,14) (2045,16) (2046,18) (2047,20) (2048,20) (2049,21) (2050,22)
};
\end{axis}
\end{tikzpicture}}\label{fig:ee_fb_investment_b}
	\caption{Investment patterns under the Electric Engagement (EE) and Falling Behind (FB) pathways over the 2030--2050 planning horizon: (a) annual effective CAPEX present values of transmission and storage investment decisions; (b) cumulative numbers of selected transmission and storage investments.}
	\label{fig:ee_fb_investment}
\end{figure}
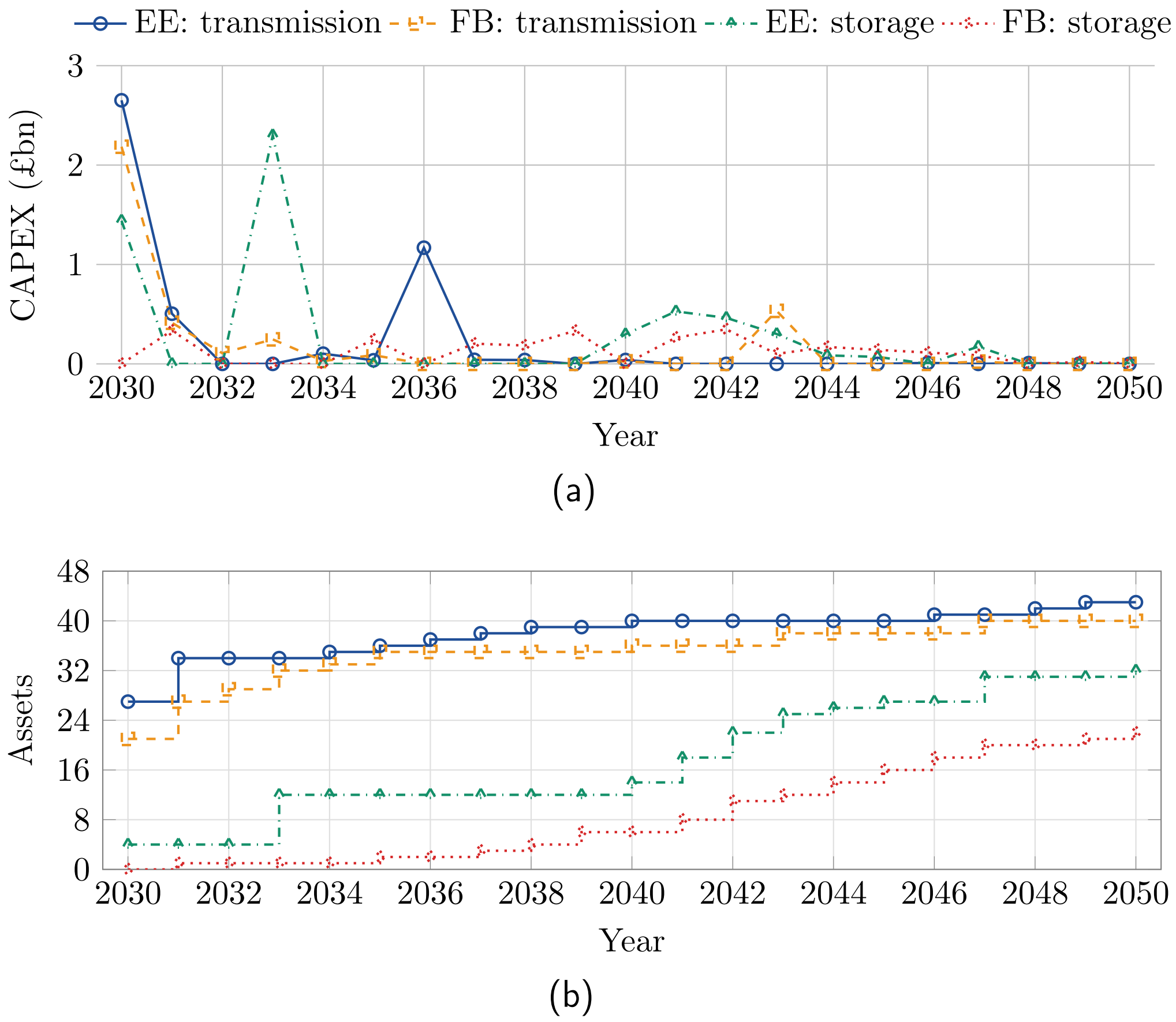

Figure~\ref{fig:ee_fb_investment} shows that transmission investment is concentrated early in the planning horizon under both pathways. EE selects 27 transmission reinforcements in 2030 and reaches 34 by 2031, compared with 21 and 27 under FB. The EGL1-2735 investment decision occurs in 2036 under EE and 2043 under FB, illustrating the earlier selection of an important reinforcement under EE. Individual
projects do not all follow the same timing pattern, and some shared
reinforcements are selected earlier under FB.

Storage deployment differs substantially between the pathways. EE selects
four storage facilities in 2030 and eight additional facilities in 2033,
reaching 12 facilities by that year. FB has no new storage investment in 2030
and selects its first facility in 2031. Subsequent storage investment extends
to 2050 in both cases, with final totals of 32 facilities under EE and
22 under FB. Thus, the pathways affect both the timing and the final scale
of storage deployment.

\subsubsection{Spatial Distribution of Expansion Decisions}
\begin{figure}[H]
	\centering
	
	\subfloat[Electric Engagement (EE).]
	{
		\includegraphics[width=0.47\columnwidth]
		{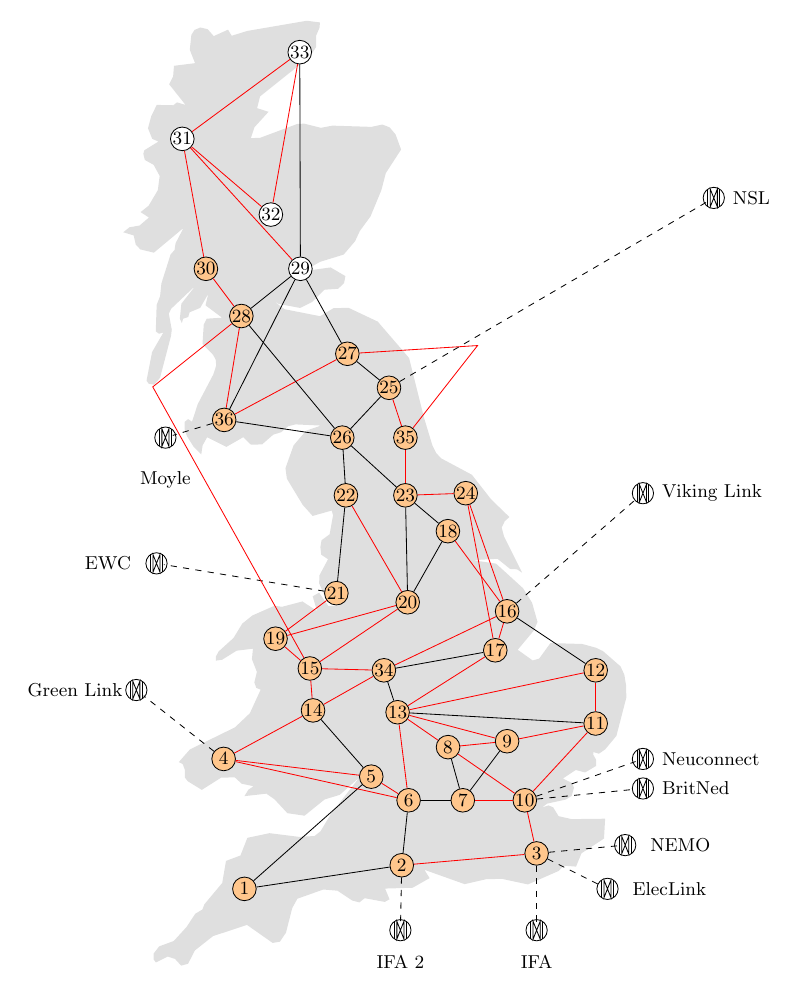}
		\label{fig:ee_network}
	}
	\hfill
	\subfloat[Falling Behind (FB).]
	{
		\includegraphics[width=0.47\columnwidth]
		{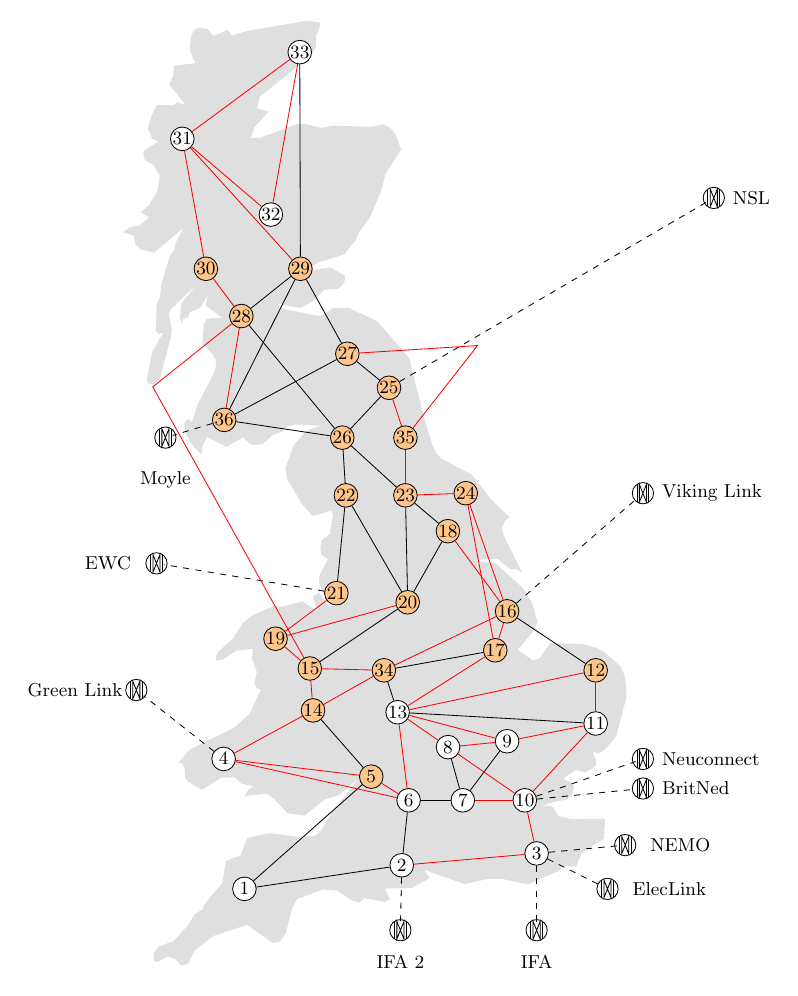}
		\label{fig:fb_network}
	}
	
	\caption{Spatial distribution of transmission reinforcement and storage
		deployment under the two FES pathways: (a) Electric Engagement (EE)
		and (b) Falling Behind (FB), aggregated over the 2030--2050 planning
		horizon. }
	\label{fig:ee_fb_network}
\end{figure}

Figure~\ref{fig:ee_fb_network} and Table~\ref{tab:case_ef_comparison} show that EE selects 43 transmission reinforcements with 66.750 GW of additional capacity, whereas FB selects 40 reinforcements with 62.650 GW. EE therefore selects three additional reinforcements and 4.100 GW, or 6.54\%, more transmission capacity. All 40 reinforcements selected under FB are also selected under EE. The three EE-specific assets are L31-1520, L43-2022, and L54-2736.

EE selects 32 storage facilities with 38.4 GWh of capacity, compared with
22 facilities and 26.4 GWh under FB. The EE storage capacity is therefore
45.45\% higher. The two portfolios share 21 storage sites; EE selects
11 additional sites not selected under FB, while FB selects one site,
S29, that is not selected under EE. The pathways therefore affect storage
location as well as capacity and timing.

\subsubsection{Scenario Implications}
The EE pathway requires both more transmission capacity and more storage
capacity than FB. Relative to FB, the increase is 6.54\% for transmission
capacity and 45.45\% for storage energy capacity. In proportional terms,
the difference in final storage deployment is therefore larger than the
difference in transmission expansion for these two pathways. This comparison
does not isolate the effects of individual scenario assumptions, but it
demonstrates that the resulting planning differences extend to both asset types.

The 40 transmission reinforcements selected under both pathways indicate a
substantial common investment portfolio, while the additional EE-specific
reinforcements and differences in storage deployment illustrate pathway-dependent
requirements. Selection under both pathways identifies common candidates;
it does not by itself establish robustness to all possible future conditions.

\subsection{Scenario-Specific and Hourly Analysis of Load Shedding}
\label{subsec:scenario_hourly_load_shedding}

The residual load shedding in Cases A and B warrants a scenario-specific
examination, particularly in view of the effectively zero load shedding
obtained in the EE-only and FB-only cases. Inspection of the nodal results
shows that all material load shedding in Cases A and B occurs under the
HT pathway at zone15 during 2045--2047. No material
load shedding occurs under EE, FB, or HE in either stochastic solution.
Thus, the results for Cases E and F are consistent with the EE and FB
outcomes within Cases A and B.

\begin{figure}[H]
	\subfloat[2045]{\definecolor{caseA}{RGB}{31,78,151}
\definecolor{caseB}{RGB}{237,148,29}

\begin{tikzpicture}
\begin{axis}[
    ybar, bar width=6pt,
    width=\textwidth, height=3cm,
    title={2045}, title style={font=\normalsize},
    symbolic x coords={0,1,2,3,4,5,6,7,8,9,10,11,12,13,14,15,16,17,18,19,20,21,22,23},
    xtick=data,xticklabels={0,1,...,23},
    xlabel={Hour}, ylabel={LS (GW)},
    ymin=0, ymax=2700, ytick={0,1000,...,3000},yticklabels={0,1,...,3},
    tick label style={font=\small}, label style={font=\small},
    ymajorgrids, xmajorgrids, grid style={dashed, gray!40},
    axis line style={draw=none}, tick style={draw=none},
    legend style={at={(0.5,1.12)}, anchor=south, legend columns=2, draw=none, font=\small},
    enlarge x limits=0.02,
]
\addplot[fill=caseA, draw=none,single ybar legend] coordinates {
(0,0) (1,0) (2,0) (3,0) (4,0) (5,0) (6,0) (7,0) (8,0) (9,0) (10,0)
(11,0) (12,0) (13,0) (14,0) (15,0) (16,0) (17,13.675447) (18,208.475275)
(19,826.205405) (20,1808.686403) (21,191.119297) (22,0) (23,0)
};
\addplot[fill=caseB, draw=none,single ybar legend] coordinates {
(0,0) (1,0) (2,0) (3,0) (4,0) (5,0) (6,0) (7,0) (8,0) (9,0) (10,0)
(11,0) (12,0) (13,0) (14,0) (15,0) (16,0) (17,50.495974) (18,1213.660431)
(19,601.662214) (20,451.223911) (21,731.119297) (22,0) (23,0)
};
\legend{Case A, Case B}
\end{axis}
\end{tikzpicture}}\\
	\subfloat[2046]{\definecolor{caseA}{RGB}{31,78,151}
\definecolor{caseB}{RGB}{237,148,29}

\begin{tikzpicture}
\begin{axis}[
    ybar, bar width=6pt,
    width=\textwidth, height=3cm,
    title={2046}, title style={font=\normalsize},
    symbolic x coords={0,1,2,3,4,5,6,7,8,9,10,11,12,13,14,15,16,17,18,19,20,21,22,23},
    xtick=data,xticklabels={0,1,...,23},
    xlabel={Hour}, ylabel={LS (GW)},
    ymin=0, ymax=2700, ytick={0,1000,...,3000},yticklabels={0,1,...,3},
    tick label style={font=\small}, label style={font=\small},
    ymajorgrids, xmajorgrids, grid style={dashed, gray!40},
    axis line style={draw=none}, tick style={draw=none},
    legend style={at={(0.5,1.12)}, anchor=south, legend columns=2, draw=none, font=\small},
    enlarge x limits=0.02,
]
\addplot[fill=caseA, draw=none,single ybar legend] coordinates {
(0,0) (1,0) (2,0) (3,0) (4,0) (5,0) (6,0) (7,0) (8,0) (9,0) (10,0)
(11,0) (12,0) (13,0) (14,0) (15,0) (16,0) (17,0) (18,602.595767)
(19,1270.083503) (20,2327.477371) (21,850.905430) (22,0) (23,0)
};
\addplot[fill=caseB, draw=none,single ybar legend] coordinates {
(0,0) (1,0) (2,0) (3,0) (4,0) (5,0) (6,0) (7,0) (8,0) (9,0) (10,0)
(11,0) (12,0) (13,0) (14,0) (15,0) (16,0) (17,114.737469) (18,2474.051525)
(19,310.636127) (20,1992.665541) (21,158.971408) (22,0) (23,0)
};
\legend{Case A, Case B}
\end{axis}
\end{tikzpicture}}\\
	\subfloat[2047]{\definecolor{caseA}{RGB}{31,78,151}
\definecolor{caseB}{RGB}{237,148,29}

\begin{tikzpicture}
\begin{axis}[
    ybar, bar width=6pt,
    width=\textwidth, height=3cm,
    title={2047}, title style={font=\normalsize},
    symbolic x coords={0,1,2,3,4,5,6,7,8,9,10,11,12,13,14,15,16,17,18,19,20,21,22,23},
    xtick=data,xticklabels={0,1,...,23},
    xlabel={Hour}, ylabel={LS (GW)},
    ymin=0, ymax=2700, ytick={0,1000,...,3000},yticklabels={0,1,...,3},
    tick label style={font=\small}, label style={font=\small},
    ymajorgrids, xmajorgrids, grid style={dashed, gray!40},
    axis line style={draw=none}, tick style={draw=none},
    legend style={at={(0.5,1.12)}, anchor=south, legend columns=2, draw=none, font=\small},
    enlarge x limits=0.02,
]
\addplot[fill=caseA, draw=none,single ybar legend] coordinates {
(0,0) (1,0) (2,0) (3,0) (4,0) (5,0) (6,0) (7,0) (8,0) (9,0) (10,0)
(11,0) (12,0) (13,0) (14,0) (15,0) (16,0) (17,0) (18,1287.173910)
(19,0) (20,1190.444760) (21,410.489961) (22,0) (23,0)
};
\addplot[fill=caseB, draw=none,single ybar legend] coordinates {
(0,0) (1,0) (2,0) (3,0) (4,0) (5,0) (6,0) (7,0) (8,0) (9,0) (10,0)
(11,0) (12,0) (13,0) (14,0) (15,0) (16,0) (17,0) (18,1557.784419)
(19,968.456094) (20,0) (21,361.868118) (22,0) (23,0)
};
\legend{Case A, Case B}
\end{axis}
\end{tikzpicture}}
	\caption{Hourly load shedding in Cases A and B under the Holistic
		Transition pathway at zone15. The upper, middle, and lower panels
		show the representative-day profiles for 2045, 2046, and 2047,
		respectively. Values are scenario-conditional load shedding in GW
		before probability weighting. The horizontal axis follows the
		model's hour indices 0--23, and all panels use the same vertical scale. }
	\label{fig:case_ab_hourly_load_shedding}
\end{figure}
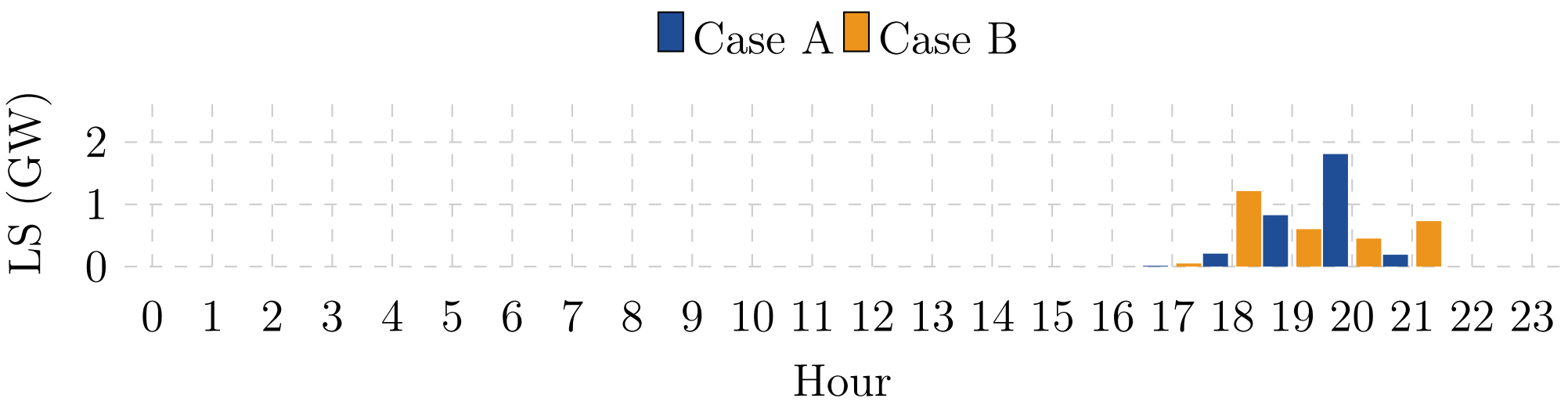
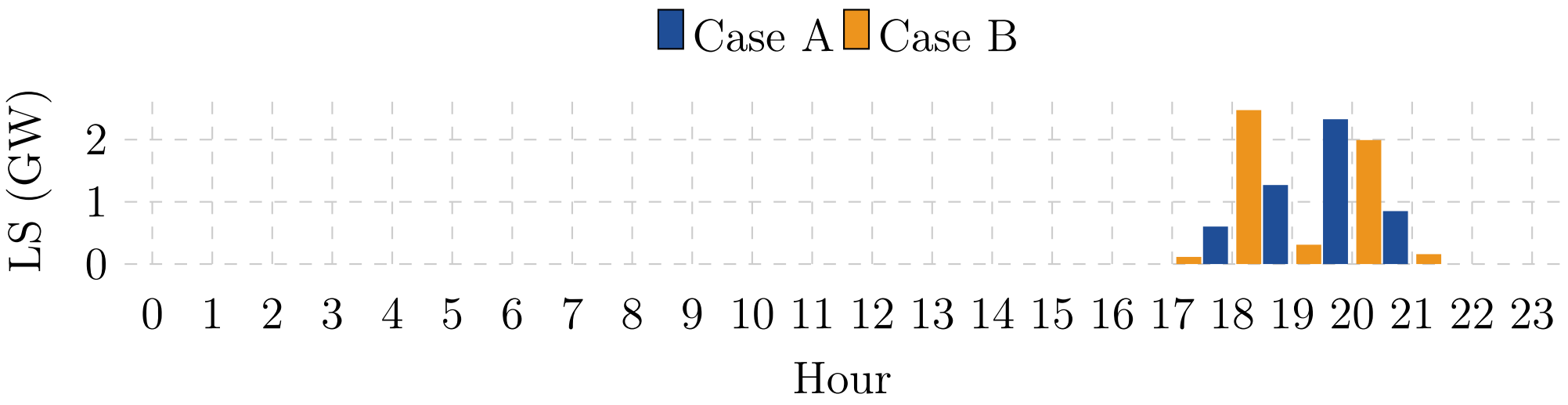
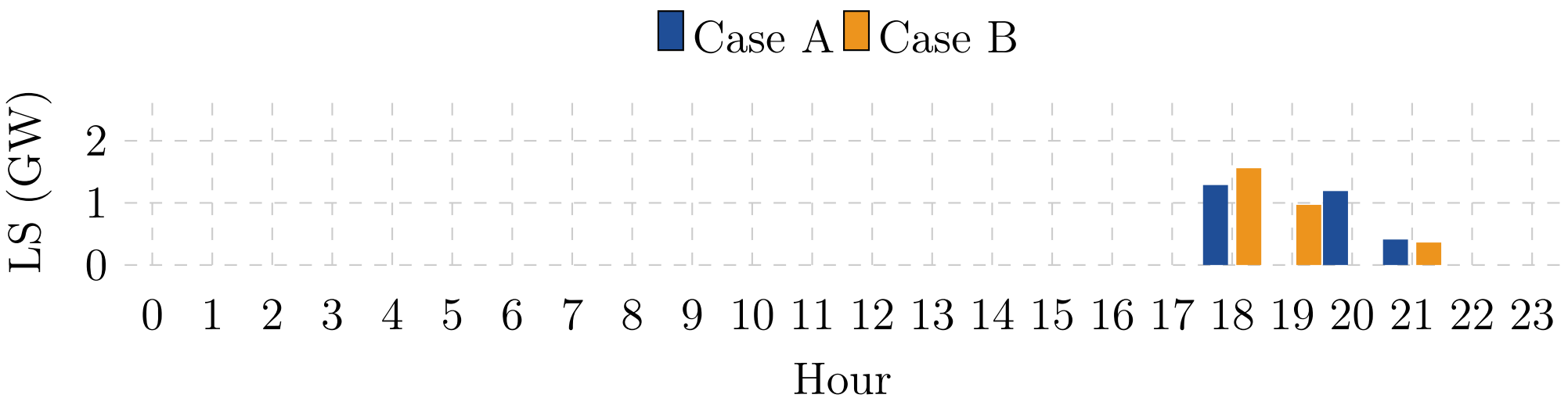

Figure~\ref{fig:case_ab_hourly_load_shedding} presents the hourly load shedding
in the representative day for each affected year. The bars show
scenario-conditional load shedding in GW, without multiplying by the HT
probability of 0.25. All 24 modeled hours are included, and grouped bars
allow the two cases to be compared without overlapping curves. Material
load shedding is concentrated within hours 17--21, although the affected
hours and the distribution of unserved energy differ between the cases.

In 2045, the maximum hourly load shedding is 1,808.69 MW at hour 20 in
Case A and 1,213.66 MW at hour 18 in Case B. In 2046, the corresponding
peaks are 2,327.48 MW at hour 20 and 2,474.05 MW at hour 18. In 2047,
both peaks occur at hour 18, reaching 1,287.17 MW and 1,557.78 MW,
respectively. Consequently, removing the transmission preparation period
does not consistently reduce the maximum hourly shortfall, even though
Case B has a lower total expected system cost.

Despite these hourly differences, the representative-day unserved energy
is equal within numerical precision between Cases A and B: approximately
3,048.16 MWh in 2045, 5,051.06 MWh in 2046, and 2,888.11 MWh in 2047.
Applying the model's annualisation factor of 365 and the HT probability
of 0.25 gives annual expected load shedding of approximately 0.278,
0.461, and 0.264 TWh. The cumulative expected value is therefore
1.003 TWh, corresponding to approximately 4.010 TWh conditional on HT.
These annual values are derived from representative-day profiles rather
than a chronological simulation of all hours in each year.

The absence of material load shedding under EE and FB does not imply
that the other pathways must also be free of shortfalls. Although EE
and FB represent contrasting energy-transition trajectories, they do not
necessarily provide upper and lower bounds on nodal supply adequacy.
Higher electricity demand can be accompanied by different generation
availability and spatial supply patterns. Consequently, neither aggregate
demand nor the overall degree of electrification is sufficient to rank
the difficulty of serving every node in every hour. The observed HT
shortfall identifies a specific location and time window requiring
attention, rather than establishing that HT is uniformly more demanding
than EE across the system.

The distinction between stochastic and single-pathway planning must also
be interpreted carefully. Cases A and B select a common investment
portfolio for all four pathways, whereas Cases E and F optimize
investments separately for EE and FB. A common portfolio may involve
trade-offs between scenario-specific requirements. However, this general
possibility does not establish the cause of the observed load shedding:
EE and FB are already served without material shortfalls by the common
portfolios in Cases A and B. The comparison therefore provides no direct
evidence that stochastic planning sacrifices EE or FB adequacy. An
HT-only investment optimisation would be needed to test whether a
portfolio tailored to HT could reduce the residual shortfall.

The Case A--B comparison further separates investment timing from annual
adequacy outcomes. Both cases select the same final transmission and
storage assets, although some investment decision years differ. Their
matching annual unserved-energy totals show that removing preparation
periods is insufficient to eliminate the HT shortfall in these solutions.
Nevertheless, the hourly profiles differ substantially: the largest
shortfall in 2046 occurs at hour 20 in Case A and hour 18 in Case B.
Annual expected load shedding alone therefore does not capture changes
in the timing or magnitude of individual hourly interruptions. These
differences describe the returned solutions and should not be interpreted
as proof that preparation periods uniquely determine the hourly allocation
of load shedding; alternative near-optimal operating schedules may also
contribute.

At the operational level, the observed shortfall means that demand at
zone15 exceeds the supply delivered by local generation, net network
inflows, and net storage discharge in the selected dispatch. Its
concentration in a limited set of hours is consistent with a local and
time-dependent adequacy issue, but does not identify which constraint
is decisive. Establishing the mechanism requires examining adjacent-line
flows and limits, local generation headroom, and storage state of charge
and charging/discharging limits in the affected hours. Minimizing HT
load shedding with the investment portfolio fixed would help distinguish
an unavoidable shortfall under that portfolio from the allocation chosen
by the cost-minimising dispatch. Allowing HT-specific investment would
then help assess the role of the common investment portfolio and the
available expansion options. Until these checks are performed, the
shortfall should not be attributed specifically to saturated transmission
lines, depleted storage, or an economic decision to accept unserved energy.

\section{Discussion}
This study demonstrates that combining future energy scenarios with explicit transmission network models allows for the translation of macro-level energy transition pathways into specific long-term infrastructure investment requirements. By integrating year-by-year supply and demand data from FES 2025, the GB36 network model (updated with boundary capability information from NESO ETYS), and the Western HVDC Link and EGL1 projects, the study jointly optimizes transmission expansion and energy storage deployment over the 2030–2050 planning horizon. This framework not only compares the total investment scale across different pathways but also identifies variations in investment timing and asset portfolios, thereby supporting the analysis of long-term grid and flexibility requirements within the context of the energy transition.

The preparation period for transmission projects influences investment timing and system costs, though this impact does not necessarily manifest as changes to the final asset portfolio. In both Case A and Case B, the same set of 45 transmission expansion projects and 32 energy storage facilities, adding 69.150 GW and 38.4 GWh of capacity, respectively, were ultimately selected; however, the decision years differed for 19 transmission investments and 8 energy storage investments. When the transmission preparation period was removed, Case B’s total expected system cost was approximately 0.96\% lower than that of Case A, even though the combined investment costs for transmission and energy storage were slightly higher. Consequently, this study concludes that the preparation period alters investment scheduling and the balance between investment and operating costs; it cannot be simply generalized that ignoring the preparation period leads to an underestimation of investment costs or inevitably changes final siting decisions. In long-term planning, differing construction and commissioning schedules can yield distinct economic outcomes, even when the final scale and location of expansions remain identical.

The joint planning of power transmission and energy storage highlights the complementary relationship between spatial transmission capacity and temporal balancing capability. Compared to Case C, which permits only transmission expansion, Case A sees a 2.40\% reduction in new transmission capacity, a 35.99\% decrease in total expected system costs, and a 97.02\% reduction in cumulative expected load shedding. This indicates that the primary benefit of the new results lies not in a substantial reduction of transmission infrastructure construction, but in a significant improvement in supply adequacy. Specifically, the total cost savings in Case A stem largely from reduced penalties for load shedding; excluding these penalties, the cost stands at £568.12 billion, slightly higher than the £566.62 billion recorded in Case C. Consequently, the additional investment and operating expenditures associated with energy storage should be evaluated alongside the value of the unserved energy avoided. Conversely, Case D, which allows only energy storage expansion, results in 244.836 TWh of cumulative expected load shedding despite the deployment of 31.2 GWh of storage capacity. This demonstrates that, given the network configuration, candidate investment options, and operational conditions of this study, the selected energy storage portfolio could not compensate for the supply inadequacy resulting from the absence of transmission expansion.

Scenario-level load-shedding results further demonstrate that the aggressiveness of the energy transition cannot serve as a direct basis for ranking power supply adequacy. Substantial load shedding in Cases A and B was confined entirely to Zone 15 under the HT pathway, occurring during specific evening hours on representative days between 2045 and 2047; notably, the same stochastic programming scheme resulted in no substantial load shedding under the EE, FB, or HE pathways. These findings align with the results from single-scenario optimisations (EE and FB), which showed virtually no load shedding. High levels of electrification can coincide with varying generation mixes and spatial supply structures; thus, aggregate demand alone is insufficient to determine that a specific pathway is harder to satisfy across all nodes and time periods. Eliminating the transmission lead time failed to resolve the aforementioned HT supply deficit, although differences in the hourly distribution and peak magnitude of load shedding persisted between Cases A and B. This indicates that the total annual expected load shedding does not fully capture the characteristics of hourly supply interruptions. However, current results are insufficient to definitively attribute this deficit to transmission congestion, energy storage depletion, or economic trade-offs within the optimisation; further analysis of power flows, generation margins, storage states, and system constraints during the relevant periods is required.

A comparison between the EE and FB pathways reveals that the future energy trajectory influences both the scale and timing of investments in transmission and energy storage. The EE pathway entails 66.750 GW of new transmission capacity, 6.54\% higher than the 62.650 GW under FB, and 38.4 GWh of energy storage capacity, 45.45\% higher than the 26.4 GWh under FB. Consequently, the relative difference in final energy storage capacity between the two pathways exceeds that of the new transmission capacity; thus, the impact of pathway variations cannot be attributed primarily to transmission requirements while assuming the scale of energy storage remains largely unchanged. As demand, generation mix, and other assumptions vary across the different pathways, these results reflect planning differences arising from the complete scenario configurations rather than isolating the causal effects of electrification or renewable energy growth.

A comparison of investment timelines also reveals that the EE scenario involves a larger scale of initial transmission expansion, with investment decisions for certain key projects made earlier. For instance, the investment decision for EGL1 is scheduled for 2036 under EE, compared to 2043 under FB. However, not all shared projects are slated for earlier construction under EE; thus, a slower energy transition should be understood as altering the investment requirements and sequencing of various projects, rather than simply delaying the entire investment plan across the board. In terms of the asset portfolio, the 40 transmission expansion projects selected under FB are all included in the EE investment set, with EE selecting three additional projects. Projects selected in both scenarios represent shared investment needs, whereas path-specific projects require ongoing assessment in light of future supply and demand fluctuations. Nevertheless, selection in both scenarios does not equate to proven robustness across all possible future outcomes.

This study employs a simplified GB36 network model aimed at identifying long-term strategic investment patterns rather than providing definitive construction plans for specific engineering projects. Where an ETYS boundary corresponds to multiple GB36 branches, boundary capacity is allocated using a simplification method, and the spatial growth of future demand and generation is mapped based on existing nodal proportions. Furthermore, the study utilizes representative 24-hour operating profiles for each planning year; annual load shedding volumes are derived by annualizing results from these representative days, meaning the model cannot fully capture extreme events or cross-day energy storage requirements associated with continuous year-round operation. Additionally, the DCOPF model does not fully account for AC network characteristics or security constraints. Consequently, the findings are best suited for comparing the relative impacts of planning strategies, energy pathways, and assumptions regarding lead times. Future research could incorporate multiple representative days or comprehensive time-series operational data, along with more granular spatial supply-demand variations and network security constraints; it could also employ fixed investment-operation analysis and single scenario investment optimisation specifically targeting HT—to further identify the precise mechanisms driving residual load shedding.

\section{Conclusion}

This paper proposes a multi-year joint planning framework for transmission and energy storage that accounts for future energy uncertainties. It integrates FES 2025 pathways, regional transmission capacity constraints calibrated against NESO/ETYS data, and an updated GB36 network model into an annual intertemporal optimisation framework spanning 2030-2050. The framework simultaneously considers transmission expansion, energy storage investment and operation, and transmission lead times associated with line lengths, thereby enabling a comparison of investment scale, location, and timing across different planning strategies and energy pathways.

The results indicate that the joint planning of transmission and energy storage combines spatial transmission capabilities with temporal regulation capabilities, thereby improving system supply adequacy. Compared to Case C, which permits only transmission expansion, Case A achieves a reduction in total expected system costs of approximately £325.00 billion (35.99\%); new transmission capacity decreases from 70.850 GW to 69.150 GW (a 2.40\% reduction); and cumulative expected unserved energy drops from 33.654 TWh to 1.003 TWh (a reduction of approximately 97.02\%). The total cost savings stem primarily from reduced penalties for unserved energy, as joint planning does not entirely eliminate the supply deficit. Although Case D, which permits only energy storage expansion, deploys 31.2 GWh of storage, it still results in approximately 244.836 TWh of cumulative expected unserved energy; this demonstrates that, under the network and candidate investment conditions of this study, expanding energy storage alone cannot compensate for the supply capacity shortfall caused by the absence of transmission expansion.

Analysis of the transmission lead time reveals that even if the final asset portfolio remains unchanged, the timing of investments and system costs are still affected. Although Case A and Case B ultimately selected the same 45 transmission expansion projects and 32 energy storage facilities, the years in which certain investment decisions were made differed. When the lead time constraint was removed, Case B’s total expected system cost was approximately £5.52 billion (or 0.96\%) lower than that of Case A. Therefore, explicitly accounting for lead times in multi-year planning facilitates a more accurate assessment of investment schedules and their economic implications. In both cases, residual load shedding was concentrated in Zone 15 under the HT pathway during specific periods between 2045 and 2047; the annual totals were largely consistent, indicating that removing the lead time constraint alone was insufficient to eliminate this supply deficit. These results also demonstrate that the extent of the energy transition does not solely determine supply adequacy at the nodal and temporal levels.

Different energy pathways simultaneously influence investment requirements for power transmission and energy storage. The EE pathway entails 66.750 GW of new transmission capacity and 38.4 GWh of energy storage, whereas the FB pathway involves 62.650 GW and 26.4 GWh, respectively. Compared to the FB pathway, the EE pathway requires 6.54\% more new transmission capacity and 45.45\% more energy storage capacity, indicating a more significant relative difference in the final scale of energy storage deployment between the two. The EE pathway also involves a larger initial scale of transmission infrastructure construction, with investment decisions for certain key projects made earlier; however, the timing varies across different projects. Consequently, the impact of energy pathways on long-term planning should be evaluated by considering asset types, investment locations, and decision years in conjunction.

This study provides a unified analytical framework for comparing investment strategies regarding transmission and energy storage in the future GB power system. Given the limitations associated with reduced-network representations, representative daily profiles, and DCOPF assumptions, the findings are primarily applicable to comparisons within long-term strategic planning. Future work could incorporate more comprehensive time-series operational data, more granular spatial variations in supply and demand, network security constraints, and project-specific timelines; furthermore, by diagnosing operational constraints and optimizing investments specifically for HT scenarios, the specific mechanisms driving residual load shedding could be further identified.

%
%
%

\bibliographystyle{elsarticle-num}
\bibliography{citation.bib}
\end{document}